\documentclass[12pt]{amsart}
\usepackage{amsmath,amssymb, euscript}
\usepackage{epsfig, graphicx}

\newtheorem{thm}{Theorem}[section]
\newtheorem{lem}[thm]{Lemma}
\newtheorem{prop}[thm]{Proposition}
\theoremstyle{definition}
\newtheorem{defn}[thm]{Definition}
\newtheorem{exmp}[thm]{Example}
\newtheorem{rem}[thm]{Remark}

\newcommand{\blackboard}[1]{\ensuremath{\mathbb{#1}}}

\newcommand{\smallcaps}[1]{\ensuremath{\textsc{#1}}}
\newcommand{\N}{\blackboard{N}}

\newcommand{\Z}{\blackboard{Z}}
\newcommand{\Zodd}{\blackboard{Z}_{\textrm{odd}}}

\newcommand{\X}{\blackboard{X}}
\newcommand{\A}{\blackboard{A}}
\newcommand{\F}{\ensuremath{\mathbf{F}}}
\newcommand{\free}{\smallcaps{Free}}
\newcommand{\core}{\smallcaps{Core}}

\begin{document}

\address{Azer Akhmedov, Department of Mathematics,
North Dakota State University,
Fargo, ND, 58108, USA}
\email{azer.akhmedov@ndsu.edu}

\begin{center} {\bf Non-amenability of R.Thompson's group \F} \end{center}

\vspace{0.6cm}

\begin{center} {\bf Azer Akhmedov} \end{center}

\vspace{1cm}

ABSTRACT: {\Small We present new metric criteria for non-amenability and discuss applications. The main application of the results of this paper is the proof of non-amenability of R.Thompson's group F.}

\vspace{1cm}

 \begin{center} {\Large Part One: Non-Amenability Criterion} \end{center}

    \vspace{0.7cm}

 \section{General non-amenability theorems}
 
 \bigskip 
 
 Amenable groups have been introduced by John von Neumann in 1929 in connection with Banach-Tarski Paradox, although earlier, Banach himself had understood that, for example, the group $\Z $ is amenable. J. von Neumann's original definition states that a countable discrete group is amenable iff it admits an additive invariant probability measure. There are many equivalent definitions of an amenable group, and the equivalences of these definitions are often respectable theorems. 
 
  \medskip 
  
  We will be using the following definition due to Følner (See\cite{F}) which, in its own turn, has many equivalent
versions

\begin{defn} Let $\Gamma $ be a finitely generated group. $\Gamma $ is called
amenable if for every $\epsilon  > 0$ and finite subsets $K, S \subseteq \Gamma $ , there exists a
finite subset $F \subseteq \Gamma $ such that $S \subseteq F$ and $\frac{|FK\backslash F|}{|F|} < \epsilon $.
\end{defn}

\medskip 

The set $F$ is called an $(\epsilon , K)$-Følner set. Very often one uses a loose term  ``Følner set”, and very often one assumes $K$ is fixed to be the symmetrized generating set. If $1\in K$, then the set $\{x\in  F | xK \subseteq F\}$
will be called {\em the interior of $F$}, and will be denoted as $Int_KF$. The
set $F\backslash Int_KF$ is called {\em the boundary of $F$} and will be denoted by $\partial _KF$.

    \bigskip

    In this section, we will state our main theorem which is the non-amenability criteria that we will be discussing in this paper. First, we introduce the notion of height function:
    
    \begin{defn} \label{defn:height} [Height function] Let $G$ be a group.  A function $h : G \rightarrow \mathbb {N}\cup \{0\}$ is called a height function on $G$, if the following conditions hold:
		
		(i) $h(xy) \leq h(x) + h(y)$ for all $x, y\in G$; 
		
		(ii) $h(x) = h(x^{-1})$ for all $x\in G$;
		
		(iii) $h(1) = 0$ where $1\in G$ denotes the identity element of $G$.
   
     \end{defn}

\medskip

  A good example of a height function $h(x)$ is the function $|x|$ representing the length of the element in the Cayley metric. 
   
    \bigskip

    \begin{thm} \label{prop:main} Let $\Gamma $ be a finitely generated group, $\xi , \eta \in \Gamma $, $\pi :\Gamma \rightarrow \Gamma /[\Gamma , \Gamma ]$ be the abelianization epimorphism. Let also $h:\Gamma \rightarrow \mathbb{N}\cup \{0\}$ be a height function, $h(\xi ) = h(\xi ^{-1}) = 1, h(\eta ) = h(\eta ^{-1}) = 0$. Assume that the following conditions are satisfied:

    $(A) \ \pi (\xi ), \pi (\eta )$ generate a subgroup isomorphic to $\Z ^2$,

   $(D)$ there exists an odd integer $p\in \Zodd $ such that for all $g\in \Gamma $,

  (i) for at least one $\delta \in \{0,1\} $, the equality
      $$h(g\eta^{\delta }\xi ^{-p}) = h(g)+ p \ (1) $$ is satisfied;

\medskip

  (ii) if for some $\delta _0\in \{0,1\}$ the equality $(1)$ is not
satisfied then $h(g\eta^{\delta _0}\xi ^{-p}) = h(g)-p$ and $h(g\eta^{\delta _0}\xi ^{p}) = h(g)+ p$, and for all $i\in \{1, 2\}, j\in [1,i], \epsilon \in \{-1, 1\}, k\in \mathbb{N}$, the inequality $$h(g\eta^{\delta _0-\epsilon i}\xi ^{-p}\eta ^{\epsilon j}\xi ^p) \geq h(g\eta^{\delta _0-\epsilon i}\xi ^{-p}) + p$$ holds;

	\medskip
	
	(iii) for all $a, b \in \{-2, -1, 1, 2\}$,  the inequality  $$h(\xi ^{-p}\eta ^{a}\xi ^{p}\eta ^{b}\xi ^{-p})\geq 3p$$ 
	holds.

	iv) for all $N\in \mathbb{N},  \epsilon \in \{-1,1\}$, $1\leq k\leq N, j_0\in \{0,1,2\}$ and $i_1, \ldots , i_k, j_1, \ldots , j_k\in \{-2, -1, 1, 2\}$, the equality  $$h(u\xi ^{-p}) =  h(u)+p$$  holds where 
	$u = \eta ^{j_0}\xi ^{-p}\eta ^{i_1}\xi ^{p}\eta ^{j_1}\ldots \xi ^{-p}\eta ^{i_k}\xi ^{p}\eta ^{j_k}$ and for all $q\in \{1, \ldots , k\}, s_q = j_0+(i_1+j_1) + \ldots + (i_q+j_q), r_q = s_q - j_q$, we have  $s_q, r_q\in \{0, 1, 2\}, r_q\neq 1 - \epsilon $ and $s_k = 1 - \epsilon$.

     Then $\Gamma $ is not amenable.
     \end{thm}

     \begin{rem} \label{rem:after} By passing to a subgroup if necessary, we may and will assume that $\xi $ and $\eta $ generate $\Gamma $. Condition D-(ii) implies that for all $g\in \Gamma $ and $\epsilon \in \{-1, 1\}$, we have $h(g\xi ^{\epsilon p})\in \{h(g) - p, h(g) + p\}$.  
			\end{rem}

    Theorem \ref{prop:main} applies to R.Thompson's group $\F $ to establish it's non-amenability. We discuss this application in later sections.

    \medskip
    
  \begin{defn}[Shifts of height functions] Let $G$ be a group and $h:G \rightarrow \mathbb{N}\cup \{0\}$ be a height function. Then for any $g\in G$ we can define a function $h_g:G\rightarrow \mathbb{N}\cup \{0\}$ by letting $h_g(x) = h(g^{-1}x)$ for all $x\in G$.
  \end{defn}
  
	\medskip
	
  \begin{rem} In general, $h_g:G\rightarrow \mathbb{N}\cup \{0\}$ is not a height function because $h_g(1)\neq 0$. But it can be viewed as a shifting of the function $h:G\rightarrow \mathbb{N}\cup \{0\}$ such that it measures the height not with respect to $1\in G$ but with respect to $g\in G$. Then condition $(D)$ can be re-written in terms of the shifts of the original height function as well.
  \end{rem}
  		
	\medskip
  
  \begin{rem} Notice that the equality $h(x\eta^{\delta }\xi ^{-p}) = h(x) + p$ does not hold for all $x\in \Gamma $ i.e. we cannot force this condition globally, simply because the height of an element is always non-negative. In some situations, it is useful to consider the height function of type $h:\Gamma \rightarrow \Z $, i.e. the image of the function is $\Z$ instead of $\mathbb{N}\cup \{0\}$. For example, let $G = Sol(2,d), d\geq 2$ denotes the free solvable group of derived length $d\geq 1$ on two generators $a,b$. Then the ``height function" $h(w(a,b)) = \sigma _b$ is naturally interesting, where $\sigma _b$ denotes the sum of exponents of $b$, in the word $w(a,b)$. For this height function the inequality $h(x\eta^{\delta }\xi ^{-p}) = h(x) + p$ is indeed satisfied, for all $x\in G, \ \delta \in [-2, 2]$ where we set $\eta = a, \xi = b^{-1}$. 
\end{rem}
 \medskip
 
 \begin{rem} Roughly speaking, for a finite subset $F\subset \Gamma $, the condition $h_g(x\xi ^{-p}) = h_g(x)+ p$ for all $x\in F$  means that in going from $x$ to $x\xi ^{-p}$ the height jumps up. In other words, the height jumps up along the $\xi ^{-p}$ shifts from the right. In groups, this happens at the expense of the height jumping down along the $\xi ^{p}$ shifts from the right for many $x$'s in (some neighborhood of) $F$, i.e. when we go from $x$ to $x\xi ^p$. However, the property of height jumping up locally along $\xi^{-p}$ shifts seems too strong in some examples of groups. Conditions (i)-(iv) of $(D)$ offer a very interesting substitute: on one hand it says that, along the $\xi ^{-p}$ shifts, the height indeed jumps up for all but at most one $x$ on a given horizontal segment of bounded length (condition $(D)$-(i)), on the other hand, for such ``an unsuccessful" element $x$, there are many elements in a certain neighborhood of $x$ such that the height actually jumps up along the $\xi ^p$ shift [condition $(D)$-(ii)]. Moreover, one can deduce that there are many elements related to $x$ in a special way where the height jumps along both the $\xi ^{-p}$ shift and the $\xi ^{-p}$ shift. Thus existence of an ``unsuccessful" element is compensated by the existence of ``super-successful" elements. 
 \end{rem}

   \begin{rem} The fact that we use arbitrary height function makes the claims of the theorems not only very strong but also provides great flexibility in applications. For example, very often very little is known about the Cayley metric of the group, so one can work with the most convenient height function instead. 
   \end{rem}
   
	We will need the following
	
	\begin{defn}[Horizontal and vertical lines] Let $H$ be a cyclic subgroup generated by $\eta $, and $H'$ be a cyclic subgroup generated by $\xi $. A subset $gH\subset \Gamma $ will be called a {\em horizontal line (passing through $g\in \Gamma $)}, and a subset of the form $gH'$ will be called a {\em vertical line (passing through $g\in \Gamma $)}.
If $x, y$ belong to the same horizontal line $L\subset \Gamma $, then we say $x$ is on the left (right) of $y$ if $x = y\eta ^n$ where $n < 0$ ($n > 0$). 
\end{defn}
 
	\begin{rem} Because of conditions $h(\eta ) = h(\eta ^{-1}) = 0$ and because of subadditivity of the height function $h$, if $x, y\in \Gamma $ are on the same horizontal line, then $h(x) = h(y)$, and even more generally, $h_g(x) = h_g(y)$, for all $g\in \Gamma $. So height is constant on a fixed horizontal line.
	\end{rem}
           
\vspace{1cm}

\section{Generalized binary trees}

\bigskip

We will need some notions about {\em binary trees} and  {\em generalized binary trees} of groups. Let $F$ be a finite subset of a finitely generated group $\Gamma $. For us, a binary tree is a tree such that all vertices have valence 3 or 1 and one of the vertices of valence 3 is marked and called a root.

\bigskip

  \begin{defn} [Binary Trees] A binary tree $T = (V, E)$ of $F$ is a finite binary tree such that  $V\subseteq F$. A root vertex of $T$ will be denoted by $r(T)$. Vertices of valence 3 are called internal vertices and vertices of valence 1 are called end vertices. The sets of internal and end vertices of $T$ are denoted by $Int(T)$ and $End(T)$ respectively.
  \end{defn}

\medskip

  \begin{defn}(see Figure 1) A generalized binary tree $T = (V, E)$ of $F$ is a finite tree satisfying the following conditions:

\medskip

(i)	All vertices of $T$ have valence  3  or 1. Vertices of valence 3 are called {\em internal vertices}, and vertices of valence 1 are called {\em end vertices}.

\medskip

(ii)	All vertices of $T$ either consist of triples (i.e. subsets of cardinality 3) or single  elements of $F$. If a vertex has valence 3 then it is a triple; if it has valence 1 then it is a singleton. For two distinct vertices $u, v\in V$, their subsets, denoted by $S(u), S(v)$, are disjoint. The union of all subsets (triples or singletons) representing all vertices of $T$ will be denoted by $S(T)$

\medskip

(iii)	One of the vertices of $T$ is marked and called the root of $T$. The root always consists of a triple and has valence 3, and it is always an internal vertex.  We denote the root by $r(T)$.

\medskip

(iv)	For any finite ray $(a_0 = r(T), a_1, a_2, \ldots , a_k)$  of $T$ which starts at the root, and for any $i\in \{0, 1, 2, \ldots , k\}$ a vertex $a_i$ is called a {\em a vertex of level $i$}.

\medskip

 (v) A finite ray $(b_0, b_1, \ldots , b_k)$ of $T$ is called {\em level increasing} if $level(b_i) > level(b_{i-1})$ for all $1\leq i\leq k$.

\medskip

 (vi) One of the elements of each triple vertex is chosen and called {\em a central element}, the other two elements are called {\em side elements}.

\medskip

(vii)	If $u$ is a vertex of level $n\geq 2$ of $T$ of valence 3, and $v, w$ are two adjacent vertices of level $n+1$, then the set of all vertices which are closer to $v$ than to $u$ form a {\em branch of  $T$ beyond vertex $u$}. $T$ has two branches beyond the vertex $u$; the second branch will consist of the set of all vertices which are closer to $w$ than to $u$

\medskip

(viii) Similar to (vii), we define the branches beyond the root vertex. So the tree consists of the root and the three branches beyond the root.
\end{defn}

\bigskip

 \begin{figure}[h!]
  \includegraphics[width=3in,height=3in]{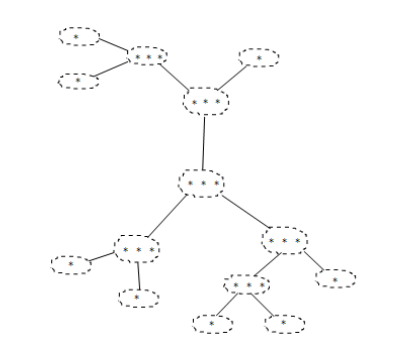}
\caption{This GBT has 6 internal vertices, and 8 end vertices; all together, it involves $3\times 6 + 8 = 26$ group elements.}
\label{labelname}
\end{figure}

\medskip

 \begin{defn}A generalized binary tree $T$ of $F\subset \Gamma $ is called trivial if it has only 4 vertices, i.e. one root vertex and 3 end vertices. $T$ is called elementary if the level of every vertex is at most 3.
 \end{defn}

\medskip

  \begin{defn}If $T = (V, E)$ is a generalized binary tree of $F\subset \Gamma $, $A\subseteq V$, then the union of all subsets (triples or singletons) which represent the vertices of $A$ will be denoted by $S(A)$.  In particular, the union of all subsets representing all vertices of $T$ will be denoted by $S(T)$.
 \end{defn}

\medskip

\begin{rem} Notice that $S(A)\subseteq F$ for all $A\subseteq V$.
\end{rem}

\medskip

 \begin{defn} The set of end vertices of a generalized binary tree $T = (V, E)$ will be denoted by $End(T)$, and the set of internal vertices will be denoted by $Int(T)$. Also, $Single(T), Central(T), Side(T)$ denote the set of all singleton vertices, central elements and side elements respectively.
 \end{defn}

\medskip

\begin{rem} By the definition of a generalized binary tree, $Single(T) = End(T)$.
\end{rem}

\medskip

 \begin{defn} For a generalized binary tree $T = (V, E)$, for all $v\in V\backslash \{root(T)\}$, $p(v)$ denotes the vertex which is adjacent to $v$ such that $level(p(v)) = level (v) - 1$; and for all $v\in V\backslash End(T)$, $n(v)$ denotes the set of vertices $v'$ which are adjacent to $v$ such that $level(v') = level (v) + 1$.
 \end{defn}

\medskip

 \begin{rem} $n$ stands for {\em next}, $p$ stands for {\em previous}. If $v$ is an internal vertex then $n(v)$ consists of pair of vertices unless $v$ is a root.
 \end{rem}

 \bigskip

 We will need the following

 \bigskip

\begin{lem} \label{prop:lemma} If $T=(V,E)$ is a GBT then

 $$   |S(End(T))| \geq \frac {1}{4}|S(T)|  $$
 \end{lem}

\medskip

 {\bf Proof.}  \ The proof is by induction on $|S(T)|$.

 \medskip

 For the trivial generalized binary tree we have $|S(T)| = 6, |S(End(T)| = 3$ so the inequality is satisfied. Let $T$ be any non-trivial GBT, and $v\in End(T), \ level(v)$ is maximal, $w = p(v)$. By definition, $v$ is a singleton and $w$ is a triple vertex.

 \medskip

 Since $level(v)$ is maximal, $n(w)$ consists of pair of singleton vertices. Let $u$ denotes the other singleton  vertex in $n(w)$. Let also $w = (a,b,c)$. We denote $w'= (a)$.

   \medskip

 By deleting $u$ and $v$ from $T$, and replacing $w$ with $w'$ we obtain anew GBT $T'$.
Then $|S(T')| = |S(T)| - 4, |S(End(T'))| = |S(End(T))| - 1 \ (\star )$. By inductive hypothesis, we have  $|S(End(T'))| \geq \frac {1}{4}|S(T')|$ which, by $(\star )$, immediately implies  $|S(End(T))| \geq \frac {1}{4}|S(T)|$ \ $\square $

    \vspace{1cm}

     {\bf Quasi-GBTs:} We need even more general objects than GBTs, namely quasi-GBTs. A quasi-GBT is a somewhat degenerate form of a GBT. The major difference is that internal vertices of odd level are allowed to be pairs (instead of triples). 

     \medskip
     
      \begin{defn} \label{defn:quasi-GBT} A quasi-GBT $T = (V, E)$ of $F$ is a finite tree satisfying the following conditions:

\medskip

(i)	All vertices of $T$ have valence  $i\in \{1, 2, 3\}$. Vertices of valence $i\in \{2, 3\}$ are called {\em internal vertices}, and vertices of valence 1 are called {\em end vertices}.

\medskip

(ii)	All vertices of $T$ consist of $k$-tuples of elements of $F$ where $k\in \{1, 2, 3\}$. If a vertex has valence $i$ then it is an $i$-tuple. For two distinct vertices $u, v\in V$, their subsets, denoted by $S(u), S(v)$, are disjoint. The union of all subsets (triples, pairs or singletons) representing all vertices of $T$ will be denoted by $S(T)$

\medskip

(iii)	One of the vertices of $T$ is marked and called the root of $T$. The root is always an internal vertex, and contains three elements.  We denote the root by $r(T)$.

\medskip

(iv)	For any finite ray $(a_0 = r(T), a_1, a_2, \ldots , a_k)$  of $T$ which starts at the root, and for any $i\in \{0, 1, 2, \ldots , k\}$ a vertex $a_i$ is called a {\em a vertex of level $i$}.

\medskip

 (v) A finite ray $(b_0, b_1, \ldots , b_k)$ of $T$ is called {\em level increasing} if $level(b_i) > level(b_{i-1})$ for all $1\leq i\leq k$.

\medskip

 (vi) One of the elements of each internal vertex is chosen and called {\em a central element}, the elements of the vertex other than central element are called {\em side elements}.

\medskip

 (vii) vertices of even level are triples.
 
 \medskip

(viii)	If $u$ is a vertex of level $n\geq 2$ of $T$ of valence $i\in \{2, 3\}$, and $v_1, \ldots , v_{i-1}$ are the adjacent vertices of level $n+1$, then the set of all vertices which are closer to $v_i$ than to $u$ form a {\em branch of  $T$ beyond vertex $u$}. $T$ has $(i-1)$ branches beyond the vertex $u$.

  \medskip
  
  (ix) Similar to (vii), we define the branches beyond the root vertex.
\end{defn}

 \medskip
 
  \begin{defn} [$\eta $-normal quasi-GBT] Let $\Gamma $ be a finitely generated group satisfying condition $(A)$. A quasi-GBT $T$ of $F\subset \Gamma $ is called  $\eta $-{\em normal} if for every vertex $v$ of $T$, $S(T)\subset xH$ for some $x\in \Gamma $.
  \end{defn}

\vspace{0.8cm}

 {\bf Labeled quasi-GBTs:} We will introduce a bit more structure on quasi-GBTs. 

\medskip

 Let $\Gamma $ be a finitely generated group satisfying condition $(A)$, $F\subset \Gamma $ be a finite subset of $\Gamma , \ F_0\subseteq F$ is partitioned into 2-element subsets $\{x, y\}$ such that $y\in \{x\xi , x\xi ^{-1}\}$. Thus every element in $F_0$ has a $\xi $-partner  which we denote by $N_{\xi }(x)$ and since $\xi $ is fixed, we will drop it and denote by $N(x)$. By definition, $N(N(x)) = x, \forall x\in F_0$.

 \medskip

  Let $T = (V, E)$ be a quasi-GBT of $F_0$. Assume that  $e\in E(T)$ is an edge connecting $v, w\in V(T)$ such that $level (w) = level (v) + 1$. Assume also that the following conditions are satisfied:

  \medskip

  (L1) there exists $a\in S(v), b\in S(w)$ such that $a = N(b)$.

  (L2) if $w$ is an internal vertex then $b$ is the central element of $w$ and $a$ is a side element of $v$.

\medskip

  Then we label the edge $e$ by the element $a^{-1}b$. We also will denote $a = start(e), b = end(e)$.
	
	\medskip
	
	\begin{defn} \label{defn:labeledquasi-GBT} A quasi-GBT $T = (V, E)$ of $F_0\subset \Gamma $ is called {\em labeled} if every edge $e = (v, w)\in E$ with $p(w) = v$ satisfies conditions (L1) and (L2) and labeled as described above.  
\end{defn}

\bigskip
	
	Now, let $r = (b_1, b_2, \ldots , b_k)$ be a finite level increasing ray in $T$, i.e. $level (b_i) = level (b_{i-1}) + 1), \ \forall i\in \{2, 3, \ldots , k\}$. We will associate an element $L(r)\in \Gamma $ to $r$.

 \medskip

   Let $e_i$ be the edge connecting $b_i$ to $b_{i+1}, 1\leq i\leq k-1$ and $a_i = start(e_i), b_i = end(e_i), 1\leq i\leq k-1$. Since $T$ is labeled, all edges $e_1, e_2, \ldots , e_{k-1}$ are labeled by some elements $l(e_1), l(e_2), \ldots , l(e_{k-1})\in \Gamma $.

  \medskip

  On the other hand, for each vertex $b_i, 2\leq i\leq k-1$ we assign $l(b_i) = end(e_{i-1})^{-1}start(e_i)$.

   \medskip

   Thus all vertices $b_2, b_3, \ldots , b_{k-1}$ are labeled by some elements $l(b_2), l(b_3), \\ \ldots , l(b_{k-1})$. [However, notice that the labeling of these vertices actually depend on $r$; if $r, r'$ are finite level increasing rays passing through the vertex $v$ and diverging at $v$ then the label of $v$ with respect to $r$ will differ from the label of $v$ with respect to $r'$].

   \medskip

   Then we associate the group element $l(e_1)l(b_2)l(e_2)\ldots l(b_{k-1})l(e_{k-1})$ to $r$ which we will denote by $L(r)$.

\medskip

      \begin{rem} Notice that in a labeled quasi-GBT we associate a word $L(r)$ to every level increasing ray $r$. Notice also that the labeling structure of a generalized binary tree of $F$ depends on the choice of non-torsion element $\xi \in \Gamma $ and subset $F_0\subset F$ which can be partitioned into pairs of $\xi $-partners, and it depends on the partitioning as well.
      \end{rem}

     \bigskip
     
    At the end of this section, we would like to introduce an important structure on quasi-GBTs, namely, the order. This notion will be crucial, in the proof of Theorem \ref{prop:main} for guaranteeing that in building super-quasi-GBTs we do not get loops, so there is no obstacle in building the trees.
  
  \medskip
  
  \begin{defn}[Order] Let $T$ be a super-quasi-GBT. Let $<$ be a linear order on the set of vertices $V(T)$ which satisfies the following conditions
  
  \medskip
  
  (i) $root(T) < v$ for all $v\in V(T)\backslash \{root(T)\}$;
  
  \medskip
   
  (ii) for all $u, v, w\in V(T)$ if $u < v$ and $w$ belongs to the branch beyond $v$ then $w < u$.
  
  Then $T$ with $<$ will be called ordered.
 \end{defn}

\medskip

 \begin{rem} Notice that if $v$ is a vertex of level $n\geq 0$ of an ordered quasi-GBT $T$, and $v_1, \ldots , v_k$ are all vertices of level $n+1$ adjacent to $v$ (i.e. $n(v) = \{v_1, \ldots , v_k)$, then we have a complete freedom in defining the order on $\{v_1, \ldots , v_k\}$. However, notice that defining the order on the set $n(v)$ for each $v$ defines the order on the whole set $V(T)$   
 \end{rem}

      \vspace{1cm}

      \section{Zigzags and other intermediate notions}

\bigskip

 In this section we will assume that $\Gamma $ is a finitely generated group satisfying condition $(A)$. A good and useful example of such a group is the group $\mathbb {Z}^2$ itself.

\medskip

  \begin{defn} [Zigzags, quasi-zigzags] \label{defn:zigzag} Let $H_{odd} = \{\eta ^{2n+1}\ | n\in \mathbb{Z}\}, \ H'_{odd} = \{\xi ^{2n+1}\ | n\in \mathbb{Z}\}$. For every $x\in \Gamma $, we call the sets $xH$ and $xH'$ {\em the horizontal lines} and {\em the vertical lines}, respectively, {\em passing through $x$}. 
	
	A sequence $Z = (x_1, x_2, \ldots , x_m)$ of elements of $\Gamma $ will be called {\em a zigzag} if for all $1\leq i\leq m-2$, \  either $x_i^{-1}x_{i+1}\in H_{odd}$ or $x_i^{-1}x_{i+1}\in H'_{odd}$; moreover, $x_i^{-1}x_{i+1}\in H_{odd}\Rightarrow x_{i+1}^{-1}x_{i+2}\in H'_{odd}$ and $x_i^{-1}x_{i+1}\in H'_{odd}\Rightarrow x_{i+1}^{-1}x_{i+2}\in H_{odd}$. 
	
	A sequence $Z = (x_1, x_2, \ldots , x_m)$ of elements of $\Gamma $ will be called {\em a quasi-zigzag} if for all $1\leq i\leq m-2$, \  either $x_i^{-1}x_{i+1}\in H$ or $x_i^{-1}x_{i+1}\in H'_{odd}$; moreover, $x_i^{-1}x_{i+1}\in H\Rightarrow x_{i+1}^{-1}x_{i+2}\in H'_{odd}$ and $x_i^{-1}x_{i+1}\in H'_{odd}\Rightarrow x_{i+1}^{-1}x_{i+2}\in H$.
	
	The number $m$ will be called {\em the length of Z}.
  \end{defn}

 \medskip

   Notice that given any horizontal line $L$ and any zigzag $Z = (x_1, x_2,\ldots , x_m)$ in $\Gamma $ with distinct elements, if $m\geq 3$ and $x_1, x_m\notin L$ then we have $|Z\cap L| = 2n, n\in \mathbb{N}\cup \{0\}$. We will need more structures on a zigzag

 \medskip

  \begin{defn} [balanced zigzags] \label{defn:balancedzigzags} A (quasi-)zigzag  $Z = (x_1, x_2, \ldots , x_m)$ is called {\em balanced} if the following conditions are satisfied:
	
	\medskip
	
	(i) for all $1\leq i\leq m-3$, and  $\epsilon \in \{-1,1\}$, we have $x_i^{-1}x_{i+1} = \xi ^{\epsilon p} \Rightarrow x_{i+2}^{-1}x_{i+3}\in \xi ^{-\epsilon p}$.
	
	(ii) either $x_1^{-1}x_2 = \xi ^{-p}$ or $x_2^{-1}x_3 = \xi ^{-p}$
	\end{defn}
	
\medskip

{\bf Example:} A zigzag 

$Z = (z, z\eta ^{n_1}, z\eta ^{n_1}\xi ^{-p}, z\eta ^{n_1}\xi ^{-p}\eta ^{n_2}, z\eta ^{n_1}\xi ^{-p}\eta ^{n_2}\xi ^{p}, z\eta ^{n_1}\xi ^{-p}\eta ^{n_2}\xi ^{p}\eta ^{n_3}, \\ z\eta ^{n_1}\xi ^{-p}\eta ^{n_2}\xi ^{p}\eta ^{n_3}\xi ^{-p})$ is balanced, for all $n_1, n_2, n_3\in \Z_{odd}$.

\begin{rem} Notice that if a zigzag $Z = (x_1, \ldots , x_m)$ is balanced then the zigzag $Z^{-1} = (x_m, \ldots , x_1)$ is not necessarily balanced.
\end{rem}

 \begin{defn} Given subsets $A, B\subseteq \Gamma $, if there exists a balanced quasi-zigzag $Z = (x_1, \ldots , x_m)$ with $x_1\in A, x_n\in B$, then we say $B$ is {\em connected to} $A$ with a quasi-zigzag $Z$. 
 \end{defn}

 {\bf Balanced Segments.} We will be working with tilings of the group $\Gamma $ into segments.
 
 \medskip
 
 \begin{defn} [Segments] Let $x, y\in \Gamma $ belong to the same horizontal line $L$, and $x$ is on the left side of $y$, i.e. there exists $n\in \N$ such that $y = x\eta ^n$. Then by $seg(x,y)$ we denote the finite set of all points (elements) in between $x$ and $y$ including $x$ and $y$, i.e. $seg(x,y) = \{x\eta ^i \ | \ 0\leq i\leq n\}$. A finite subset $I$ of $L$ is called a segment if there exists $x, y\in I$ such that $I = seg(x,y)$.  
\end{defn}

\medskip

\begin{defn} [Balanced Segments] A finite segment $I$ will be called {\em a balanced segment} if $|I|$ is divisible by 6.
  \end{defn}
 
\medskip
 
  \begin{defn} [Leftmost and rightmost elements] Let $I$ be a segment. Then there exists a unique element $z_1\in I$ such that for any $z\in I$ there exists a non-negative integer $n$ such that $z = z_1\eta ^n$. The element $z_1$ will be called {\em the leftmost element} of $I$. Similarly, we define the rightmost element of $I$: there exists a unique element $z_2\in I$ such that for any $z\in I$ there exists a non-positive integer $n$ such that $z = z_2\eta ^n$; \ $z_2$ will be called {\em the rightmost element} of $I$
  \end{defn}

  \medskip
	
	\begin{defn}[Compatibility of zigzags and quasi-zigzags with tiling] Let $\{I(x)\}_{x\in X}$  be a tiling of $\Gamma $ by balanced segments. A (quasi-)zigzag $Z=(x_1, x_2, \ldots x_n)$ is called compatible with this tiling if for all $i\in \{1, \ldots , n-1\}$, if $x_i, x_{i+1}$ lie on the same horizontal line then there exists $x\in X$ such that $x_i, x_{i+1}\in I(x)$.
	\end{defn}

\medskip

  \begin{defn} [Compatibility of GBT with tiling] Let $\{I(x)\}_{x\in X}$  be a tiling of $\Gamma $ by balanced segments. A quasi-GBT $T$ of $\Gamma $ is called compatible with the tiling, if for all $v\in V(T)$, there exists $x\in X$ such that $S(v)\subset I(x)$.
  \end{defn}

\medskip

 We would like to conclude this section with introducing some notions for labeled $\eta $-normal generalized binary trees.
  
  \medskip

  \begin{defn} [starting element] Let $T$ be a labeled $\eta $-normal quasi-GBT, and $a$ be the central element of $root(T)$. The element $N(a)$ will be called {\em the starting element} of $T$ and denoted by $start(T)$.  
 \end{defn}

  \begin{defn} [special quasi-GBTs] A labeled $\eta $-normal quasi-GBT is called special if $\{a'\}$ is an end vertex where $a'$ is the starting element of $T$. 
  \end{defn}

	\medskip
	
	In the proof of Theorem \ref{prop:main} the GBTs will be constructed piece-by-piece, in other words, some GBTs will be constructed as a union of elementary pieces.  An trivial labeled $\eta $-normal GBT $T$ is called an {\em elementary piece}. If $T$ is special as a quasi-GBT then it is called {\em a special elementary piece}, otherwise we call it {\em an ordinary elementary piece}.   
	
	\medskip 
	
	 It is useful to observe that quasi-GBTs can be obtained from a GBT as follows: Let $T$ be an $\eta $-normal and labeled GBT. Let also $v_1, \ldots , v_k\in Int(T)$ be not necessarily distinct internal vertices such that  $v_i\neq root(T)$ for any $i\in \{1, \ldots , k\}$. Let $a_i\in S(v_i), 1\leq i\leq k$, such that $level(N(a_i)) = level(a_i) + 1$, $level(a_i)\in \Z_{odd}$, and let $w_i\in V(T)$ be such that $N(a_i)\in S(w_i)$. Finally, let $S_i =\{v\in V(T) \ | \ v$ belongs to some level increasing ray $r$ which starts at $w_i\}$. By deleting $\displaystyle \mathop{\sqcup }_{1\leq i\leq k} \{a_i\}\sqcup S(S_i)$ from $T$ we obtain an $\eta $-normal and labeled quasi-GBT. If $T$ is a special GBT then we obtain a special $\eta $-normal and labeled quasi-GBT.  
	
	\begin{defn} Let $T$ be a labeled $\eta $-normal quasi-GBT. We say $T$ is {\em successful} if it contains a triple vertex of odd level.
	\end{defn}
	
  \vspace{1cm}
    
  \section {Generalized binary trees in partner assigned regions}
     
  \medskip
  
   In this section, we will be assuming that $\Gamma $ is a finitely generated group satisfying conditions $(A)$ and $(D)$, $K = B_{2000p}(1)\subset \Gamma $, $\epsilon = \frac{1}{100|K|}$, $F$ is a connected $(K, \epsilon )$-Følner set, $\{I(x)\}_{x\in X}$ be a collection of pairwise disjoint balanced segments of length at most 1200 tiling $\Gamma $. Let also $X_0 = \{x\in X \ | \ I(x)\cap Int_KF \neq \emptyset \},  F_1 = \displaystyle \mathop{\cup }_{x\in X_0}I(x)$. 
   
  \medskip
  
  We will need the following notions
   
  \medskip
  
  \begin{defn}[Regions] A subset $S\subseteq F_1$ is called a region if there exists $X_0'\subset X$ such that $S = \displaystyle \mathop{\sqcup }_{x\in X_0'}I(x)$. For any subset $S'\subset F_1$ we will denote the minimal region containing $S'$ by $R(S')$. 
\end{defn}

 \medskip
 
 \begin{defn} [Partner assigned regions] A region $S\subseteq F_1$ is called {\em partner assigned} if there exists a subset $S_0\subseteq S$ and a function $n : S_0 \rightarrow \Zodd $ such that $$S\subseteq \displaystyle \mathop{\sqcup }_{x\in S_0} \{x, x\xi ^{n(x)}\} \subseteq F$$ We will denote  $x\xi ^{n(x)} = N(x)$ and $x = N(x\xi ^{n(x)})$. Also, for any subset $A\subseteq S$, we will write $N(A) = \{N(x) \ | \ x\in A\}$.
 \end{defn}
 
 \medskip
 
  Notice that the partners of elements from $F_1$ do not necessarily lie in $F_1$ but they always lie in $F$. In the proof of Theorem \ref{prop:main} we will be assigning partners of elements of $F_1$ before we start building the trees, so $F_1$ will be a partner assigned region. We will arrange the partner assignment to satisfy certain conditions to enable us to push the rays of the trees to higher levels or at least not to let them come below certain level.
        
  \medskip

  \begin{defn} [Zigzags and GBTs respecting partner assignment] Let $S\subseteq F_1$ be a partner assigned region and $Z = (x_1, x_2, \ldots , x_n)$ be a quasi-zigzag in $S$. We say $Z$ respects partner assignment of $S$ if for all $i\in \{1, \ldots , n-1\}$,  if $z_i$ and $z_{i+1}$ do not lie on the same horizontal line then $z_{i+1} = N(z_i)$.
	
	 We say a quasi-GBT $T$ is respecting the given partner assignment if it is a labeled quasi-GBT with respect to this partner assignment.
	 	
		A partner assignment induces labeling structure for a quasi-GBT $T$; from now on we will identify the notion of ``labeled quasi-GBT" with the notion of a ``quasi-GBT respecting the partner assignment". We will also identify the notion of ``balanced (quasi-)zigzag" with the notion of a ``(quasi-)zigzag respecting the partner assignment", if the partner assignment is fixed. (Notice that if a quasi-zigzag satisfies our fixed partner assignment, then it is balanced; cf. Definitions \ref{defn:zigzag} and  \ref{defn:balancedzigzags}.)   
  \end{defn}
    	
  \medskip

 \begin{defn} An element $z\in F_1$ is called {\em successful} if at least one of the following three conditions are satisfied:

 (i) $h(z\xi ^{-p}) = h(z) + p$.

 (ii)  $N(z)\notin F_1$ 
	
 (iii) $N(z) = z\xi ^p$

 Otherwise, $z$ is called {\em unsuccessful}.

 More generally, for any region $S\subseteq F_1$, we say $z\in S$ {\em is a successful element of} $S$, if either at least one of conditions (i), (iii) are satisfied, or  $N(z)\notin S$.       
  \end{defn}

	\bigskip
 
 Now we will introduce a certain partner assignments for regions in $F_1$. First, we describe a certain natural tiling and partner assignment of the group $\Z ^2$. Namely, for every $m,n\in \Z$, we let $I_{m,n} = \{(x,m)\in \Z^2 \ | \ 360n \leq x < 360(n+1)$. The partner assignment will be defined as follows: for every $(x,y)\in \Z$, $N(x,y) = (x,y+p)$ if $y$ is even, and $N(x,y) = (x,y-p)$ if $y$ is odd.
 
 \medskip
 
 \begin{rem} \label{rem:tiling} In the group $\Z ^2$ with standard generating set $\eta = (1,0), \xi = (1,0)$, let us assume the above tiling $\{I_{m,n}\}_{m,n\in \Z}$ and the partner assignment, and let $p = 1, h(x,y) = |y|, \mathrm{for \ all} \ (x,y)\in \Z^2$. Then for every quasi-zigzag $Z$ in $\Z^2$ which starts at element $g\in \Z^2$, and respects the given partner assignment, we have $h(z)\in (h(g)-1, h(g), h(g)+1)$ for all $z\in Z$.
\end{rem}
 
 \medskip
 
 \begin{rem} The above tiling and partner assignment immediately induces the pull back tiling by balanced segments of length $360$ and partner assignment through the epimorphism $\pi :\Gamma \rightarrow \Z ^2$ which satisfies the following conditions: 

 (i) for all $x\in \Gamma $ either $N(x) = x\xi ^p$ or $N(x) = x\xi ^{-p}$.

 (ii) for all $\epsilon \in \{-1,1\}$, if $x,y\in \Gamma $ are on the same horizontal line and $N(x) = x\xi ^{\epsilon p}$ then $N(y) = y\xi ^{\epsilon p}$. 

 \medskip
 
 We will fix this tiling in $\Gamma $ and the partner assignment for the rest of the paper.  Notice that any zigzag $Z = (x_1, \ldots , x_m)$ which respects this tiling and the partner assignment is necessarily balanced if $x_1^{-1}x_2 = \xi ^{-p}$ or $x_2^{-1}x_3 = \xi ^{-p}$.
 \end{rem}

 \bigskip

 We now would like to introduce more structures on the partner assigned regions as well as on quasi-GBTs in partner assigned regions. 

 \medskip

\begin{defn} [Suitable/unsuitable segments] A segment $I(x), x\in X$ is called {\em suitable} if for some (consequently, for any) $g\in I(x)$ we have $\pi (g) = (x,y)\in \Z ^2$ where $y\in \Z_{odd}$. If $I(x)$ is not suitable then it is called {\em unsuitable}.
\end{defn}

\medskip

\begin{rem} Notice that, by the definition of the partner assignment, the segment $I(x)$ is suitable iff for some (consequently, for any) $g\in I(x)$ we have $N(g) = g\xi ^{-p}$.
\end{rem} 

\medskip

\begin{defn} Let $T$ be an $\eta $-normal quasi-GBT, $v$ be a vertex of it, $J$ be a horizontal segment. We say $v$ belongs to $J$ if $S(v)\subseteq J $.
   \end{defn}

\bigskip

 {\bf {\em For the rest of the paper} } we make the following assumptions: let $\Gamma $ be a finitely generated amenable group satisfying conditions $(A)$ and $(D)$, $K = B_{2000p}(1)\subset \Gamma $, $\epsilon = \frac{1}{100|K|}$, $F$ be a connected $(K, \epsilon )$-Følner set, $\{I^k(x)\}_{x\in X, 0\leq k\leq 359}$ be a tiling of $\Gamma $ by segments of length 3; and for all $x\in X$, let $I(x)$ be a balanced segment of length $360\times 3 = 1080$ such that $I(x) = \displaystyle \mathop{\sqcup } _{0\leq k\leq 359}I^k(x)$; we assume that the tiling $\{I(x)\}_{x\in X}$ is obtained by the pullback of the tiling in Remark \ref{rem:tiling}.  Let also $X_0 = \{x\in X \ | \ I(x)\cap Int_KF \neq \emptyset \},  F_1 = \displaystyle \mathop{\cup }_{x\in X_0}I(x)$, and let $\overline{F_1} = F_1\cup N(F_1)$. All the labeled $\eta $-normal GBTs and quasi-GBTs will respect this tiling $\{I(x)\}_{x\in X}$ and the fixed partner assignment, and all the elementary pieces will respect the tiling $\{I^k(x)\}_{x\in X, 0\leq k\leq 359}$. The notion of the region will be understood with respect to this tiling, i.e. a subset $S \subseteq F_1$ will be called a region if there exists a subset $Y\subseteq X\times \{0,\ldots , 359\}$ such that $S = \displaystyle \mathop{\sqcup }_{(x,k)\in Y}I^k(x)$.
     
     \medskip
     
     For every $x\in X$, let $I(x) = \{z(x)\eta ^i \ | \ 0\leq i\leq 1079\}, z(x)\in \Gamma $.  Then we let $I^k(x) = \{z(x)\eta ^i \ | \ 3k\leq i < 3k+3\}, \ 0\leq k\leq 359, J^l(x) = \{z(x)\eta ^j \ | \ 6l \leq j < 6l+6\}, \ 0\leq l\leq 179$. The intervals $I^k(x), 0\leq k\leq 359$ will be called {\em short intervals}, and the intervals $J^l(x), 0\leq l\leq 179$ will be called {\em medium intervals}. All the labeled $\eta $-normal quasi-GBTs will respect the tiling $\{I^k(x)\}, 0\leq l\leq 359, x\in X$ unless specifically stated otherwise in which case they will respect the tiling $\{J^l(x)\}, 0\leq l\leq 179, x\in X$.

     \medskip

     We will say that a short segment $I^k(x), x\in X, 0\leq k\leq 359$ is suitable (unsuitable) if the segment $I(x)$ is suitable (unsuitable). We will write $Y_0 = \{(x,k) \ | \ x\in X_0, 0\leq k\leq 359\}$. We will assume that all the labeled GBTs and quasi-GBTs have a root consisting of a suitable short segment.
		
		\medskip
		
		We may and will assume that there is no unsuitable segment $I^k(x)\subset F_1$ such that for some $z\in I^k(x)$, $N(z)\notin F_1$. Also, we will assume that $F_1$ has at least as many suitable medium intervals as non-suitable ones. 
		
		\medskip
		
		If $u, v$ are on the same horizontal line and $u = v\eta ^n, n\in \Z$, we will write $d(u,v) = |n|$. For $g\in \Gamma , \ |g|$ will denote the length of $g$ in the left invariant Cayley metric of $\Gamma $ w.r.t. the generating set $\{\xi , \eta \}$.

 \bigskip
    
    \begin{figure}[h!]
  \includegraphics[width=3in,height=3in]{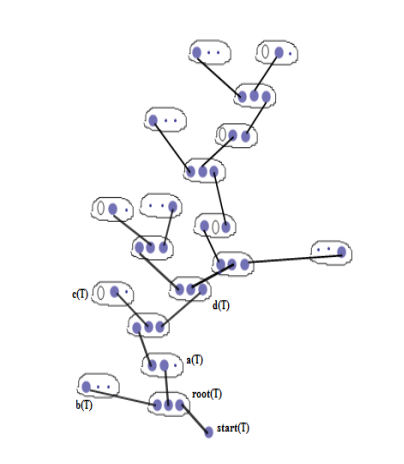}
\caption{This is a special labeled quasi-GBT with a root vertex $\mathrm{root}(T)$, starting element $\mathrm{start}(T)$, and different types of vertices $a(T),  b(t),  c(T)$ and $d(T)$. A vertex is a $k$-tuple ($k\in \{1,2,3\}$) iff it contains $k$ big black dots; an empty circle indicates that the element is omitted, and the small dot means the element is not in $F_1$.}
\label{labelname}
\end{figure}

 \bigskip

  \vspace{1cm}

  \section{Completeness}
  
  \medskip

 We make the assumptions of Section 4.
 
 \medskip
 
 In this section we introduce a key notion of completeness of quasi-GBTs. The completeness of the sequence $T_1, T_2, \ldots T_m$ simply means that we do not stop the trees unnecessarily. 
 
 \medskip

\begin{defn}[Complete sequence of quasi-GBTs.] Let $T_1, T_2, \ldots , T_n$ be a sequence
of mutually disjoint labeled $\eta $-normal quasi-GBTs of $F$. We say this sequence is pre-complete if for every $i\in \{1, \ldots , n\}$ 

 (i) the root of $T_i$ consists of a triple $I^k(y)\subset F_1, y\in X, 0\leq k\leq 359$, where $I(y)$ is a suitable segment; 

 (ii) the singleton vertices (i.e. end vertices) of $T_i$ belong to unsuitable segments; 

 (iii) if $v$ is a vertex of even level of $T_i$ then $v$ is a triple and $S(v) = I^l(z)$ where $I^l(z)\subset F_1, z\in X, 0\leq l\leq 359$;  

 \medskip

 $T_1, T_2, \ldots , T_n$ is called complete if, in addition, for all $v =\{x\}\in Single(T_i)$ with $x\in I^k(x_0)\subseteq F_1, x_0\in X, 0\leq k\leq 359$, at least one of the following conditions hold

 \medskip

 (i) there exists $j\in \{1,\ldots ,i-1\}$ such that for some $w\in V(T_j)$, $I^k(x_0)\backslash \{x\} = S(w)$;

 \medskip

 (ii) $T_i$ is a special quasi-GBT, $I(x_0)$ is a non-suitable segment, and $x$ is a starting element of $T_i$;

 \medskip

 (iii) there exists $j\in \{1, \ldots , i-1\}, x'\in I^k(x_0)\backslash \{x\}$ such that $T_j$ is a special quasi-GBT and $x'$ is the starting element of $T_j$.  
\end{defn}

\medskip

\begin{defn}[Successful complete sequence] Let $T_1, T_2, \ldots , T_n (\ast )$ be a complete sequence
of labeled $\eta $-normal quasi-GBTs of $F$ and $I$ be an unsuitable short segment in $F_1$. We say $(\ast )$ is successful at $I$ if one of the following conditions hold:

 (i) there exist distinct $z_1, z_2\in I$ such that $z_1 = start (T_{i}), z_2 = start (T_{j})$ for some $i, j\in \{1, \dots , n\}$.
 
 ii) there exist $i, j\in \{1, \dots , n\}$ such that $start (T_i)\sqcup S(v) = I$ where $v$ is a pair vertex of $T_j$.
 
 iii) there exist $i\in \{1, \dots , n\}$ such that $S(v) = I$ for some triple vertex $v$ of $T_i$.
 
\end{defn}

 Notice that if a pre-complete sequence $T_1, T_2, \ldots , T_n$ satisfies the condition $\displaystyle \mathop{\sqcup }_{1\leq i\leq m}S(T_i) \supseteq F_1$ then it is complete unless there exists an unsuitable segment $I = I^k(x) = \{z_1, z_2, z_3\}\subset F_1$ such that for all $i\in \{1,2,3\}$, there exists $j\in \{1, \dots , n\}$ such that $\{z_i\}$ is an end vertex of $T_j$, moreover, if $T_j$ is a special quasi-GBT, then $z_i$ is not a starting element of $T_j$. In this case, we will say that the sequence $(T_1, \dots , T_m)$ is {\em unsuccessful at $I$}.
 
 \medskip 
 
  It will be useful to observe that labeled quasi-GBTs can be obtained from a labeled GBTs by simply deleting some branches. Let $T$ be a labeled GBT such that all the internal vertices of $T$ are triples of the form $(z, z\eta ^{n_1}, z\eta ^{n_2})$ for some $n_1, n_2\in \Z $; moreover, all end vertices are of even level. Let $v_1, \ldots , v_k\in Int(T)$ such that $level(v_i)\in \Z_{odd}$ for all $1\leq i\leq k$. Let also $a_i\in S(v_i), 1\leq i\leq k$ such that level $(N(a_i)) = level(a_i) + 1$, and let $w_i\in V(T)$ be such that $N(a_i)\in S(w_i)$. Finally let $S_i = \{v\in V(T) \ | \ v \ \mathrm{belongs \ to \ some \ level}$ \ \ $\mathrm{increasing \ ray } \ r \ \mathrm{which \ starts \ at} \ w_i\}$. By deleting $\sqcup _{1\leq i\leq k} \{a_i\}\sqcup S(S_i)$ from $T$ we obtain {\em an ordinary labeled quasi-GBT} $R$. If $T$ is special (or respects the partner assignment, or a tiling) then we will say that $R$ is special (or respects the partner assignment, or a tiling). Notice that the definition of labeled quasi-GBT as in the above construction agrees with the definitions \ref{defn:quasi-GBT} (quasi-GBT) and  \ref{defn:labeledquasi-GBT} (labeled quasi-GBT).

 \medskip

 The following proposition will be very useful in the proof of Theorem \ref{prop:main}.  

\medskip

\begin{prop} \label{prop:complete} Let $T_1, T_2, \ldots , T_m$ be a complete sequence of labeled $\eta $-normal quasi-GBTs of $F$ such that $\displaystyle \mathop{\sqcup }_{1\leq i\leq m}S(T_i) \supset F_1$. Assume that for every ball $B_{20p}(g), g\in F_1$ there exists $i\in \{1, \ldots , m\}$ such that $T_i$ is successful and $S(T_i)\cap B_{20p}(g) \neq \emptyset $. Then $|\partial _KF| \geq \frac {1}{100|B_{50p}|}|F|$. $\square $
\end{prop} 

\medskip

  {\bf Proof.} Let us recall that for every $x\in X$, $I(x) = \{z(x)\eta ^i \ | \ 0\leq i\leq 1079\}$, where $z(x)\in \Gamma $ is the leftmost element of $I(x)$, and $I^k(x) = \{z(x)\eta ^i \ | \ 3k\leq i < 3k+1\}$. We also have $J^l(x) = I^{2l}(x)\sqcup I^{2l+1}(x)$ for all $x\in X, 0\leq l\leq 179$. Notice that $\{J^l(x)\}_{x\in X, 0\leq l\leq 179}$ forms a tiling of $\Gamma $; also every $J^l(x)$ has length 6 and contains exactly two short intervals; we will call these short intervals {\em neighbors}.
  
  \medskip
  
  Let $Z_1 = \{2l+1 \ | \ 0\leq l\leq 179\}, Z_2 = \{2l \ | \ 0\leq l\leq 179\} , F_1^{-} = \sqcup _{x\in X_0, k\in Z_2}I^k(x),  F_1^{+} = \sqcup _{x\in X_0, k\in Z_1}I^k(x)$. Notice that $F_1 = F_1^{-}\sqcup F_1^{+}$. Notice also that, since the quasi-GBTs $T_1, \ldots , T_m$ respect the tiling $\displaystyle \mathop{\sqcup }_{x\in X_0, 0\leq k\leq 359}I^k(x)$, for all $i\in \{1,\ldots , m\}$, either $S(T_i)\cap F_1^{-} = \emptyset $ or $S(T_i)\cap F_1^{-} = \emptyset $. We may assume that if $v\in V(T_i), 1\leq i\leq m$, and $S(v)\backslash F_1\neq \emptyset $ then $v$ is an end vertex; moreover, $F_1$ contains at least as many suitable segments as non-suitable ones. We will also assume that $F_1^{-}$ contains at least as many unsuitable short intervals at which $T_1, T_2, \ldots , T_m$ successful as $F_1^{+}$ (and we will ignore the successfulness in $F_1^{+}$). We may furthermore assume that for some $1\leq m' < m$, $S(T_i)\cap F_1^{-} = \emptyset $ for all $1\leq i\leq m'$, and  $S(T_i)\cap F_1^{+} = \emptyset $ for all $m'+1\leq i\leq m$.

 \bigskip

   As agreed in Section 4, labeled $\eta $-normal quasi-GBTs $T_1, \dots , T_m$ respect the tiling $\{I^k(x)\}_{x\in X, 0\leq l\leq 359}$. Notice that for such a quasi-GBT $T$, if $I_1, I_2$ are two neighboring short intervals, then for all $1\leq i\leq m$, either $S(T_i)\cap I_1 = \emptyset $ or $S(T_i)\cap I_2 = \emptyset $. In the proof below, we will construct labeled $\eta $-normal quasi-GBTs will respect just the tiling $\{J^l(x)\}_{x\in X, 0\leq l\leq 179}$. Recall that by the agreement the root and more generally, all vertices of even level of any labeled $\eta $-normal quasi-GBT consists of a suitable short interval.
  
  \medskip
  
  Let $J_i = I_i^{-}\sqcup I_i^{+}, 1\leq i\leq N$ be all medium intervals in $F_1$ such that $I_i^{-}, I_i^{+}$ are the neighboring short intervals and at $I_i^{-}$ the sequence $T_1, \ldots , T_m$ is successful. We will call the intervals $J_i, 1\leq i\leq N$ {\em successful}.
  
  \medskip
  
  {\em The main idea} is to construct GBTs $\Omega = (\overline{T}_1, \ldots , \overline{T}_M)$ of $F$ such that the following conditions hold:
  
  $(c1):$ every suitable short interval in $F_1$ forms an internal vertex of one of the GBTs from $\Omega $;
   
   \medskip
   
  $(c2):$ every unsuitable medium interval in $F_1$ contains an internal vertex of one of the GBTs from $\Omega $;
   
   \medskip
   
  $(c3)$ every successful medium interval contains either two internal vertices or one internal vertex and one starting element of GBTs from $\Omega $.
   
  \medskip
  
   Let now $T$ be a labeled $\eta $-normal quasi-GBT and $R$ be a labeled quasi-GBT such that $S(R)\cap S(T) = \emptyset $.  We introduce the following types of operations:
  
  \medskip
  
  {\bf {\em breaking up a quasi-GBT:}} Let $v$ be a pair vertex of $T, \  S(v) = \{a,b\}, level(a) = level(N(a))+1, level(b) = level(N(b))-1$. Let also $w$ be a vertex of $T$ such that $N(b)\in S(w)$. Then $w$ is a triple vertex and is adjacent to two vertices $w_1, w_2$. Let $V_i =\{v\in V(T) \ | \ v$ belongs to some level increasing ray $r$ which starts at $w_i\}$ for all $i\in \{1, 2\}$. 
  $S_i =\{v\in V(T) \ | \ v$ belongs to some level increasing ray $r$ which starts at $w_i\}$. Then $S:= \{b\}\sqcup S(w)\sqcup S(w_1)\sqcup S(w_2)\sqcup S(V_1)\sqcup S(V_2)$ forms a labeled  $\eta $-normal quasi-GBT, while $S(T)\backslash S$ forms another labeled $\eta $-normal quasi-GBT. Hence, $T$ can be broken into two labeled $\eta $-normal quasi-GBTs. If we break $T$ at all of its $k$ pair vertices (notice that $k\leq 3$) then we obtain $k+1$ $\eta $-normal quasi-GBTs. We will call these {\em $\eta $-normal quasi-GBTs obtained as a result of break up}.  
  
  \medskip

   {\bf {\em gluing of type 1:}} Let $v$ be as above, $R$ be special, $\sigma _1 = \{c_1\}$ be its starting vertex such that $v$ and $\sigma _1$ belong either to the same or neighboring short intervals. Then $S(T)\sqcup S(R)$ forms a new labeled quasi-GBT $R'$ where $root(R') = root(T)$, and the pair vertex $v$ is replaced with a triple vertex.  
   
   \medskip

   {\bf {\em gluing of type 2:}} Let $R$ be special, $\sigma _1 = \{c_1\}$ be its starting vertex and  $\sigma _2 = \{c_2\}$ be the end vertex of $T$ such that $c_1, c_2$ belong to the same short interval. Then the union $S(T)\sqcup S(R)$ forms a new labeled quasi-GBT $R'$ such that $root(R') = root(T)$.
   
   \medskip
   
   {\bf {\em collecting special quasi-GBTs:}} Let $R_1, R_2, R_3$ be disjoint special labeled quasi-GBTs, with starting elements $s_1, s_2, s_3$ respectively. Then the set $S(R_1)\sqcup S(R_2)\sqcup S(R_3)$ forms a GBT $R'$ with $root(R') = (s_1, s_2, s_3)$. All internal vertices of $R'$ consists of a triple on some horizontal line except perhaps the root.  
    
   \medskip
    
  Now, we are ready to start the proof. Inductively, for any $i\in \{1,\ldots , m'\}$ we will associate a finite set of labeled quasi-GBTs $C_i = (S^{(i)}_1, \ldots , S^{(i)}_{n_i})$ to the collection $T_1, \ldots , T_i$ as follows:
 
 \medskip
 
  {\em Base:} We let $n_1 = 1, S^{(1)}_1 = T_1$. 

 \medskip
 
  {\em Step:} Assume that, for some $i\in \{1, \ldots , m'-1\}$, $C_i = \{S^{(i)}_1, \ldots , S^{(i)}_{n_i}\}$ is the set of labeled quasi-GBTs associated to $T_1, \ldots , T_i$ such that $\displaystyle \mathop{\sqcup }_{1\leq q\leq n_i}S(S^{(i)}_q) =  \displaystyle \mathop{\sqcup }_{1\leq q\leq i}S(T_q)$.   
  
  \medskip
  
  Let $\Sigma = \{\sigma _1, \ldots , \sigma _k\}$ be a maximal subset of $End(T_{i+1})$ such that the following two conditions are satisfied:
  
  (i) if $1\leq q_1 < q_2 \leq k$ then $\sigma _{q_1}$ and $\sigma _{q_2}$ do not belong to the same short interval;
    
  (ii) for all $q\in \{1, \ldots , k\}$, $\sigma _q$ belongs to $F_1$.
  
  \medskip
  
  If $\Sigma = \emptyset $ then we let $C_{i+1} = C_i \sqcup T_{i+1}$. Otherwise, we have $\sigma _1: = \{c_1\}$ belongs to $F_1^{-}$. Let also $\sigma _1$ belongs to the short interval $I$. Then we have two cases:
  
  \medskip
  
  {\em Case 1:} $I\backslash \displaystyle \mathop{\sqcup }_{1\leq q\leq n_i}S_q^{(i)}$ contains an element besides $c_1$.
  
  \medskip
  
  Then, by completeness, we have two sub-cases:
  
  a) \ There exists $j\in \{1, \ldots , i\}$ such that the starting element $s$ of $T_j$ belong to the same short interval $I$ as $\sigma _1$. 
  
  \medskip
  
  Then we perform gluing of type 2, between $\sigma _1$ and $s$.
   
  \medskip
  
  b) \ $I\cap \displaystyle \mathop{\sqcup }_{1\leq q\leq n_i}S_q^{(i)} = \emptyset $ and $c_1$ is the starting element of $T_{i+1}$.
  
  \medskip
  
  In this case we do not perform any operation and go to $\sigma _2$.
  
  \medskip
  
  {\em Case 2:} $I\backslash  \displaystyle \mathop{\sqcup }_{1\leq q\leq n_i}S_q^{(i)}$ contains no element besides $c_1$.
  
  \medskip
  
   Then we do not perform any operation and by going to $\sigma _2$, apply the process to $\sigma _2$, and so on.
  
  \medskip
  
  Once we are done with all the list $\{\sigma _1, \ldots, \sigma _k\}$ we obtain a new set $C_{i+1}$. We denote the final collection $C_{m'} = \{S^{(m')}_1, \ldots , S^{(m')}_{n_{m'}}\}$ by $D = \{S_1, \ldots , S_n\}$.
  
  \medskip
  
  Notice that as a result of the construction, the collection $D = \{S_1, \ldots , S_n\}$ of labeled quasi-GBTs satisfies the following conditions:
  
  \medskip
  
  (i) $S_{q_1}\cap S_{q_2} = \emptyset $, for all $1\leq q_1 < q_2\leq n$;

  \medskip
  
  (ii) $\displaystyle \mathop{\sqcup }_{1\leq q\leq n}S(S_q) = \displaystyle \mathop{\sqcup }_{1\leq q\leq m'}S(T_q)$;
  
  \medskip
   
  (iii) any short interval in $F_1^{-}$ contains either a pair vertex or a triple vertex of one of the quasi-GBTs from $D$, or it consists of three starting elements of some three quasi-GBTs from $D$; 
  
  \bigskip
  
  Now, inductively on $j\in \{0, 1, \ldots , m-m'\}$, we associate a finite set of labeled quasi-GBTs $D_j = \{R^{(j)}_1, \ldots , R^{(j)}_{k_j}\}$ to the collection $T_{m'}, T_{m'+1}, \ldots , T_{m'+j}$ as follows: 
  
  \medskip
  
  {\em Base:} $D_0 = D$.
  
  \medskip
  
  {\em Step:} For $j\in \{1, \ldots , m-m'\}$, we perform break up operations at all pair vertices of $T_{m'+j}$, and let $\Phi _1, \ldots , \Phi _r$ be a set of all quasi-GBTs and $\Psi  = \{s_1, \ldots , s_r\}$ be the set of starting elements obtained as a result of the break up. Let also $K_0\subseteq \{1, \ldots ,r\}$ be the unique maximal subset such that for all $t\in K_0$ if $s_t$ belongs to the short interval $I$ then $I\cap \displaystyle \mathop{\sqcup }_{1\leq q\leq k_{j-1}}S(R^{(j-1)}_{q}) = \emptyset $.
  
  \medskip
  
  For all $t\in K_0$, let $s_t$ belongs to the short interval $I_t$ (recall that $I_t\subset F_1^{+}$). Let $I_t'$ be the neighboring short interval in $F_1^{-}$.
  
  \medskip
  
  Then we have one of the following two cases: 
  
  \medskip
  
  {\em Case 1:} $I_t'$ contains a pair vertex of some labeled quasi-GBT $S$ from $D = D_{0}$. 
  
  \medskip
  
  Then we perform gluing of type 1 between $S$ and $\Phi _t$.
  
  \medskip
  
  {\em Case 2:} $I_t'$ contains a triple vertex of some labeled quasi-GBT $S$ from $D = D_{0}$. 
  
  \medskip
  
  Then we do not perform any operation and let $\Phi _t$ be one of the quasi-GBTs in the collection $R_{j}$ (so some of the quasi-GBTs $(\Phi _t)_{t\in K_0}$ gets glued to one of the earlier quasi-GBTs from $D_{0}$ and some remain intact).
  
  \medskip
  
   Then, on the final collection $D_{m-m'} = \{R^{(m-m')}_1, \ldots , R^{(m-m')}_{k_{m-m'}}\}$, we perform the following operation: if any short interval $I\subset F_1^{-}$ consists of three starting elements of some three labeled quasi-GBTs then we collect these labeled quasi-GBTs into one GBT. We denote the resulting collection of GBTs by $\bar{D}$. Notice that, conditions (c1)-(c3) are satisfied for $\bar{D}$, i.e.: 
   
   \medskip
   
   (i) every suitable short interval in $F_1$ forms an internal vertex of one of the GBT from $\bar{D}$;
   
   \medskip
   
   (ii) every unsuitable medium interval in $F_1$ contains an internal vertex of one of the GBT from $\bar{D}$;
   
   \medskip
   
   (iii) every successful medium interval contains either two internal vertices or one internal vertex and one starting element of GBTs from $\bar{D}$.
   
   \medskip
   
   Finally, we perform the following final operations on $\bar{D}$: if $\bar{D}$ contains $z = 3z_1+z_2$ GBTs ($z_1, z_2\in \Z , z_1\geq 1, 0\leq z_2\leq 2$) $\Xi _1, \ldots , \Xi _z$ with starting elements in $F_1$ then for all $q\in \{0, \ldots , z_1-1\}$ we collect   $\Xi _{3q+1}, \Xi _{3q+2}, \Xi _{3q+2}$ into a new GBT. 
   
   \medskip
   
   Let $R_1, \ldots , R_n$ denotes the set of GBTs obtained as a result of the process. Then the following conditions hold:
   
   \medskip
   
   (i) for every suitable medium interval $J\subseteq F_1$, $J\subset \displaystyle \mathop{\sqcup }_{1\leq i\leq n}S(Int(R_i))$;
   
   \medskip
   
   (ii) for every unsuitable medium interval  $J\subseteq F_1$, $|J\cap \displaystyle \mathop{\sqcup }_{1\leq i\leq n}S(Int(R_i))|\geq 3$;
   
   \medskip
   
   (iii) for every successful unsuitable medium interval  $J\subseteq F_1$, $|J\cap  \displaystyle \mathop{\sqcup }_{1\leq i\leq n}S(Int(R_i))|\geq 4$. 
   
   \medskip
   
   Then, $|F_1\cap \displaystyle \mathop{\sqcup }_{1\leq i\leq n}S(End(R_i))| \leq \frac {1}{4}|F_1| (\ast )$, moreover, since the number of suitable medium intervals is not less than the number of unsuitable ones, we obtain that $|F_1\cap \sqcup _{1\leq i\leq n}S(Int(R_i))| \geq \frac {3}{4}|F_1| + \frac {1}{3}L$ where $L$ denotes the number of successful medium intervals.  
    
    \medskip
    
    Let $\Omega $ be a maximal subset of $F_1$ such that for any two distinct $x, y\in \Omega $, $d(x,y)\geq 70p$. Then $B_{35p}(x)\cap B_{35p}(y) = \emptyset $, moreover, $F_1\subseteq \cup _{x\in \Omega }B_{70p}(x)$. Then $|\Omega | \geq 
\frac {|F_1|}{|B_{70p}|}$.
    
     \medskip
     
    Now, by the assumption every ball of radius $20p$ centered at $F_1$ has a non-empty intersection with a successful medium interval. Then every ball of radius $35p$ centered at $F_1$ (in particular, the balls $B_{35p}(x), x\in \Omega )$  contains a successful medium interval, thus we obtain that $L \geq |\Omega | \geq \frac {|F_1|}{|B_{70p}|}$. 
    
    \medskip

    Then, by Lemma \ref{prop:lemma} we obtain that $$|\sqcup _{1\leq i\leq n}S(End(R_i))|\geq \frac {1}{3}(\frac {3}{4}|F_1|+\frac {1}{3}L) \geq \frac {1}{4}|F_1| + \frac {1}{9}L$$ 
    
    Then, by $(\ast )$, $|\sqcup _{1\leq i\leq n}S(End(R_i))\backslash F_1| \geq \frac{1}{9}L \geq \frac {|F_1|}{9|B_{70p}|}$. Hence $|\partial _KF| \geq \frac {|F_1|}{9|B_{70p}|}$ \ $\square $

  \bigskip
  
  \begin{figure}[h!]
  \includegraphics[width=3in,height=3in]{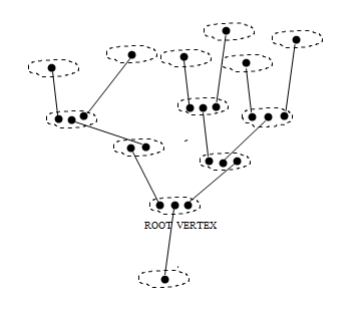}
\caption{this labeled quasi-GBT has only one pair vertex.}
\label{labelname}
\end{figure}
  
  \bigskip
  
  \begin{figure}[h!]
  \includegraphics[width=3in,height=3in]{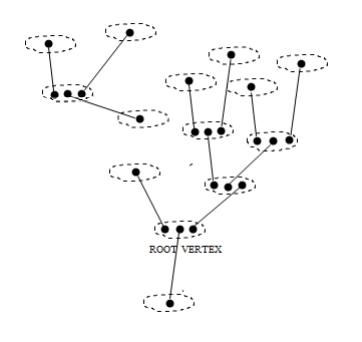}
\caption{after the break up at the pair vertex, we obtain two new labeled quasi-GBT; one of them is a trivial GBT.}
\label{labelname}
\end{figure}

 \bigskip

   The following proposition follows immediately from the proof of Proposition \ref{prop:complete}
 
 \begin{prop} \label{prop:+-} Let $T_1, T_2, \ldots , T_m (\ast )$ be a complete sequence of quasi-GBTs of $F$ such that $\displaystyle \mathop{\sqcup }_{1\leq i\leq m}S(T_i) \supset F_1$. Assume that there exist $m+n$ short unsuitable intervals $I_1, \ldots , I_m, I_{m+1}, \dots , I_{m+n}$ in $F_1$ such that $(\ast )$ is successful at $I_j$ for all $j\in \{1, \ldots , m\}$ and unsuccessful at $I_j$ for all $j\in \{m+1, \ldots , m+n\}$, moreover, $(\ast )$ is not unsuccessful at all other short unsuitable intervals of $F_1$. If $m-n \geq \frac {1}{100|B_{50p}|}|F|$, then $|\partial _KF| \geq \frac {1}{900|B_{50p}|}|F|$. \ $\square $
 \end{prop}	
    
    \vspace{1cm}
    
    \section{The Proof of Theorem \ref{prop:main}.}

  We make the assumptions of Section 4. The following notions will be needed
	
\medskip

  \begin{defn} Let $x_1, x_2\in X, A_1\subseteq I(x_1), A_2\subseteq I(x_2)$. We say $A_1, A_2$ are connected by a quasi-zigzag if for some $u_1\in A_1, u_2\in A_2$ there exists a quasi-zigzag $Z$ such that $Z$ starts at $u_1$ and ends at $u_2$. If $Z$ is also balanced, then we say $A_2$ is connected {\em to} $A_1$ by a balanced quasi-zigzag. 
		\end{defn}
		
		\medskip
	
	\begin{defn} Let $(x_1, k_1), (x_2, k_2), (x_3, k_3), (x_4, k_4)\in Y_0$ and $I_i = I^{k_i}(x_i), 1\leq i\leq 4$. We say the pairs $(I_1, I_2)$ and $(I_3, I_4)$ are connected with non-interfering zigzags if there exist zigzags $Z_1, Z_2$ in $F_1$ such that $I_1$ and $I_2$ are connected by $Z_1$, $I_3$ and $I_2$ are connected by $Z_2$, moreover, there exists no $(x,k)\in Y_0$ such that  $Z_i\cap I^{k}(x) \neq \emptyset $ for all $i\in \{1,2\}$. 
	\end{defn}
	
	\medskip

	We will also need the following notions
	
	\begin{defn} Let $S\subseteq F_1$. We write $X(S) = \{(x,k)\in X_0\times \{0, \ldots , 359\} \ | \ I^k(x)\subset S, I(x) \mathrm{ \ is \ suitable }\}, P(S) =  \displaystyle \mathop{\sqcup }_{(x,k)\in X(S)}I^k(x)$.
	\end{defn}
	
	\begin{defn}[Extremal regions] Let $g_1, g_2, \ldots , g_n\in \Gamma $ be distinct elements. Let the sub-region $S\subseteq F_1$ be such that there exist a sequence $(S_0, S_1, \ldots , S_n)$ of subsets and sequences $(h_1, \ldots , h_n)$ and $(q_1, \ldots , q_n)$ of non-negative integers such that

\medskip
 
 (i) $S_0 = P(F_1)$;

\medskip

 (ii) $h_i = max\{h_{g_i}(z) \ | \ z\in S_{i-1}\}$, for all $1\leq i\leq n$;

\medskip

 (iii) $S_i = \{z\in S_{i-1} \ | \ h_{g_i}(z) \geq h_i-q_i\}$, for all $1\leq i\leq n$,

\medskip

 (iv) $S = S_n$.

\medskip

 We will write $R(g_1, \ldots g_n; q_1, \ldots , q_n) = S\cup N(S)\cup (\overline{F_1}\backslash F_1)$, and call $R(g_1, \ldots g_n; q_1, \ldots , q_n)$ the extremal region of the sequence $(g_1, \ldots , g_n)$. 
  \end{defn}

\begin{rem} In the above definition, $S_n$ are defined inductively. Notice that the set $S_n\cup N(S_n)$ may contain a large region even if $q_1 = q_2 = \dots = q_n = 0$. In general, since $h_1, \dots , h_n$ are arbitrary non-negative integers, by taking them sufficiently big we also obtain that the entire $F_1$ (more precisely, $F_1\cup N(F_1)\cup (\overline{F_1}\backslash F_1) = F_1\cup (\overline{F_1}\backslash F_1)$) is an extremal region. This observation will be used in the sequel, but the most interesting case of an extremal region is when $(q_1, \ldots , q_n) = (0, \ldots , 0)$; this observation will be crucial in our study. A somewhat more general case of $0 \leq q_i \leq 2p$, for all $1\leq i\leq n$ is also interesting but we will not be using it in this paper. Notice also that $S_n$ is a union of suitable segments. Moreover, although $S_n$ is a region, $R(g_1, \ldots g_n; q_1, \ldots , q_n)$ may not be; we hope this (calling it an extremal region) will not cause a confusion.
\end{rem}

\begin{rem} We will sometimes denote the extremal region by $R(g_1,\ldots , g_n)$ dropping the sequence $(q_1, \ldots , q_n)$. Also, if $(q_1, \ldots , q_n) = (0, \ldots , 0)$ then we will denote the extremal region by $R_0(g_1, \ldots , g_n)$. 

 The set $R(g_1, \ldots , g_n)$ heavily depends on the order of the sequence; for example, in general, $R(g_1, g_2)\neq R(g_2, g_1)$; in fact, it is quite possible that $R(g_1, g_2)\cap R(g_2, g_1) \subseteq \bar{F_1}\backslash F_1$.
\end{rem} 

\medskip

 \begin{defn} [Minimal elements] Let $I(x), x\in X$ be a suitable segment in the extremal region $R(g_1, \ldots g_n)$, and $I\subseteq I(x)$ be a sub-segment. We say $z\in I$ is {\em a minimal element of} $I$ if $z$ is unsuccessful with respect to some $h_{g_i}, 1\leq i\leq n$, and there is no $u\in I$ and $j\in \{1, \ldots , i-1\}$ such that $u$ is unsuccessful with respect to $h_{g_j}$. We say $z$ is an {\em absolutely minimal element of} $I$, if it is minimal, and for every other minimal $u\in I$, 

 $\mathrm {min}\{j : z \ \mathrm {is \ unsuccessful \ w.r.t.} \ h_{g_j}, \ \mathrm {but} \ u \ \mathrm {is \ successful \ w.r.t.} \ h_{g_j}\} < \mathrm {min}\{j : u \ \mathrm {is \ unsuccessful \ w.r.t.} \ h_{g_j}, \ \mathrm {but} \ z \ \mathrm {is \ successful \ w.r.t.} \ h_{g_j}\}$    
\end{defn}

\bigskip
   
    Condition (i) of $(D)$ states that on a given horizontal segment of length 2 at most one element is unsuccessful with respect to the same $h_{g}$ . Then condition (ii) states that given an unsuccessful element $u$, one can relate to it an element $u'$ in a certain special way such that not only $u'$ is successful but even nicer property holds for $u'$, namely, $h(u'\xi ^{-p}\eta ^{i_1}\xi^p\eta ^{i_2}\xi ^{-p}\eta ^{i_3}\xi ^{p}) = h(u') + 4p$ for certain values of $i_1, i_2, i_3$. Thus existence of an unsuccessful element $u$ is compensated by the existence of $u'$ which is related to $u$ and is even better than successful. 
    
   \medskip
   
   We will materialize this idea in the proof. First, we need the following notions
   
   \medskip
   
   \begin{defn}[Successfully related elements] Let $u\in I^k(x)\subset F_1 $ be a suitable short interval, and $u\in I^k(x)$. We say an element $u'$ is {\em successfully related to} $u$ if either $u' = u\eta ^{-i}\xi ^{-p}\eta ^j$, where $i\in \{1, 2\}, j\in [1,i], u\eta ^{-i}\in I^k(x)$ or $u' = u\eta ^{i}\xi ^{-p}\eta ^{-j}$, where $i\in \{1, 2\}, j\in [1,i], u\eta ^i\in I^k(x)$. \end{defn}
     
  \medskip

  \begin{defn}[Successfully related special elementary pieces] Let $(g_1, \ldots , g_N)$ be a finite sequence of elements of $\Gamma $, and $u\in I = I^k(x) \subset F_1, x\in X_0, 0\leq k\leq 359$, be such that $u$ is an absolutely minimal element of $I$. Let $T$ a labeled $\eta $-normal special elementary piece of $F$. We say $T$ is {\em successfully related to} $u$ if $N(u)$ is a starting element of $T$ (so, in particular, it is an end vertex). In this case we also say $T$ is successfully related to $I$.
  \end{defn}
  
	\begin{defn} Let $T$ be a special elementary piece successfully related to $u\in I$ for some suitable short segment $I = I^k(x) = \{u, u_1, u_2\}, x\in X_0, 0\leq k\leq 359$. Let $A\subseteq F$. We say $T$ is successful relative to $A$, if for all $i\in \{1, 2\}$, if $N(u_i)$ belongs to an unsuitable short segment $I_i\subset F_1$, then  $A$ contains all elements of $I_i$ which are successfully related to $u$. 
	\end{defn}

 \medskip

   Let us emphasize that if a suitable short segment $I = I^k(x) = \{z_1, z_2, z_3\}\subset F_1$, where $z_1$ and $z_3$ are the leftmost  and rightmost elements of $I$ respectively, in the extremal region $R(g_1, \ldots , g_N)$ admits a minimal element but not an absolutely minimal element, then $I$ has exactly two minimal elements, namely, $z_1$ and $z_3$. Let $T$ be an ordinary elementary piece such that $root(T) = I$; we will denote the vertex of $T$ which contains $N(z_i)$ by $v_i, 1\leq i\leq 3$. Let also $S(v_i)\subseteq I_i = I^{k_i}(y_i), 1\leq i\leq 3$.

\medskip

 Then both elements of $I_2\backslash S(v_2)$ are successfully related to either $z_1$ or $z_2$; $N(z_1)\eta $ is successfully related to $z_3$ and $N(z_3)\eta ^{-1}$ is successfully related to $z_1$. These observations are useful but in the sequel we will choose the sequence $(g_1, \dots , g_n)$ dense enough that every suitable short segment will contain an absolutely minimal element.

 \bigskip

  Now, we would like to introduce a notion of partial (or linear) order 
	in the set of suitable short intervals of $F_1$ as well as in the set of all elements belonging to the suitable short intervals of $F_1$. The proof of Theorem \ref{prop:main} is based in constructing a complete sequence of elementary pieces or quasi-GBTs,  and the order we introduce will indicate where to start the next quasi-GBT in the sequence at every step.

   \medskip
   
   {\bf {\em Partial Order:}} We now introduce a strict partial order in extremal regions [we say a relation $\sim $ on a set $S$ is {\em a strict partial order} if it is  anti-symmetric, transitive, and there is no $x\in S$ such that $x\sim x$]. By prolonging the sequence of origins, this partial order can be made a linear order. 
   
   \medskip
   
    Let $R = R(g_1, \ldots g_n; q_1, \ldots , q_n)$ be an extremal region, and \\   $S(g_1, \ldots g_n; q_1, \ldots , q_n)$ be the set of all elements of $R(g_1, \ldots g_n; q_1, \ldots , q_n)$ which belong to suitable segments. Let also $Y_R = \{(x,k) \ | \  I^k(x)\subset S(g_1, \ldots g_n; q_1, \ldots , q_n)$. We will introduce a partial order $\prec $ on the set $Y_R$ and on the set $S(g_1, \ldots g_n; q_1, \ldots , q_n)$ as follows:
    
 \medskip
 
  Let $y_1, y_2\in Y_R, y_1 = (x_1, k_1), y_2 = (x_2, k_2)$. Let also $u_i\in I^{k_i}(x_i), 1\leq i\leq 2$. We say $y_2\prec y_1$ if one of the following two conditions hold:
  
  \medskip
  
  (c1) $\mathrm{min}\{i\in \{1, \ldots , n\} \ | \ h_{g_i}(u_1) > h_{g_i}(u_2)\} \leq \mathrm{min}\{i\in \{1, \ldots , n\} \ | \ h_{g_i}(u_2) > h_{g_i}(u_1)\}$, moreover, there exists $i\in \{1, \ldots , n\}$ such that $h_{g_i}(u_1)\neq h_{g_i}(u_2)$;
  
  \medskip
  
  (c2) $h_{g_i}(u_1) = h_{g_i}(u_2)$ for all $i\in \{1, \ldots , n\}$, but
   
  $\mathrm{min}\{i\in \{1, \ldots , n\} \ | \ I^{k_1}(x_1) \ \mathrm{contains \ an  \ unsuccessful \ element \ w.r.t.} \ h_{g_i}\} < \mathrm{min}\{i\in \{1, \ldots , n\} \ | \ I^{k_2}(x_2) \ \mathrm{contains \ an  \ unsuccessful \ element \ w.r.t.} \ h_{g_i}\}$ (so, in particular, there exists $i\in \{1, \ldots , n\}$ such that $I^{k_1}(x_1)$ contains an unsuccessful element w.r.t. $h_{g_i}$).
  
  \medskip
  
  Notice that, since the functions $h_g:\Gamma \rightarrow \mathbb{N}\cup \{0 \}$ are constant on the horizontal lines, the definition does not depend on the choices $u_1, u_2$. If $y_2\prec y_1$ then we also say $I^{k_1}(y_1)$ {\em is bigger than} $I^{k_2}(y_2)$. We also say $I^{k_1}(y_1)$ {\em is strongly bigger than} $I^{k_2}(y_2)$ if condition (c1) holds. 
  
  \medskip
  
  Let now $u, v\in S(g_1, \ldots g_n; q_1, \ldots , q_n), u\in I^{k'}(x'), v\in I^{k''}(x''), y' = (x', k'), y'' = (x'', k'')$. If $y'\prec y'' $ then we let $u\prec v$. If $y'\nprec y''$ and $y''\nprec y'$ (in particular, if $y' = y''$), we say $u\prec v$, if \\ $\mathrm{min}\{i\in \{1, \ldots , n\} \ | \ u \ \mathrm{is \ unsuccessful \ w.r.t.} \ h_{g_i} \mathrm{but} \ v \ \mathrm{is \ successful \ w.r.t.} \ h_{g_i}\} > \\ \mathrm{min}\{i\in \{1, \ldots , n\} \ | \ v \ \mathrm{is \ unsuccessful \ w.r.t.} \ h_{g_i} \ \mathrm{but} \ u \  \mathrm{is \ successful \ w.r.t.} \ h_{g_i} \}$. 
  
  \medskip
  
  Once we define the order $\prec $ on the set of suitable short segments of an extremal region (and on the set of elements which belong to suitable short segments of the extremal region) we can restrict it to any subset. Then it is clear from the condition D-(iii) that for any region $R\subseteq F_1$ (in particular, for the region $F_1$ itself) we can choose long enough sequence $g_1, \ldots , g_L$ of origins such that for the choice $(q_1, \ldots , q_L) = (0, \ldots , 0)$ the partial order $\prec $ on the set of suitable short segments of $R$ (and on the set elements which belong to suitable short segment of $R$) becomes a linear order. Indeed, in defining the partial order $\prec $, we observe that $$R_0(g_1)\supseteq R_0(g_1, g_2)\supseteq \dots \supseteq R_0(g_1, g_2, \dots , g_L).$$ Notice that $R_0(g_1, g_2, \dots , g_L)$ is never empty. If it contains only one suitable short segment of $F_1$ then this segment is the biggest in our order. But if it contains more than one such segment then we look into the question of which one of them contains unsuccessful elements with respect to $h_{g_i}$ for the least possible $i$. This way of viewing our definition of the partial order motivates more general notion of {\em extremal regions with constraints}. We will avoid defining this notion in its most general natural version, but rather restrict ourselves to a certain type (two types) which will be used in the sequel. 
  
   First, we need a the following
   
   \begin{defn} Let $Z = (x_1, \dots , x_n)$ be a zigzag in $F_1$ which respects the fixed tiling and partner assignment. We say $Z$ belongs to the class $\mathcal{R}_{left}$ ($\mathcal{R}_{right}$) if $x_i$ is not the leftmost (rightmost) element of any unsuitable short segment for all $1\leq i\leq n-1$ such that $x_i^{-1}x_{i+1} = \xi ^p$.   
   \end{defn}

   \begin{defn}[Extremal regions with constraints on the left] Let $g_1, g_2, \ldots , g_n\in \Gamma $ be distinct elements and $I_1, I_2, \dots , I_n$ be a sequence of suitable short segments of $F_1$. Let the sub-region $S\subseteq F_1$ be such that there exist a sequence $(S_0, S_1, \ldots , S_n)$ of subsets and sequences $(h_1, \ldots , h_n)$ and $(q_1, \ldots , q_n)$ of non-negative integers such that

\medskip
 
 (i) $S_0 = P(F_1)$;

\medskip

 (ii) $h_i = max\{h_{g_i}(z) \ | \ z\in S_{i-1}\}$, for all $1\leq i\leq n$;

\medskip

 (iii) for all $1\leq i\leq n$, $I_{i}\subseteq  \{z\in S_{i-1} \ | \ h_{g_i}(z) \geq h_i-q_i\}$, 

\medskip 

 (iv) $S_i = \{z\in S_{i-1} \ | \ h_{g_i}(z) \geq h_i-q_i , z \ \mathrm{is \ connected \ to } \ I_i \ \mathrm{by \ a \  zigzag \ from} \ \mathcal{R}_{left}\}$, for all $1\leq i\leq n$,

\medskip

 (v) $S = S_n$.

\medskip

 We will write $R(g_1, \ldots g_n; q_1, \ldots , q_n; I_1, \dots , I_n) = S\cup N(S)\cup (\overline{F_1}\backslash F_1)$, and call $R(g_1, \ldots g_n; q_1, \ldots , q_n; I_1, \dots , I_n)$ the extremal region of the sequence $(g_1, \ldots , g_n)$. 
 
 \medskip 
 
  By replacing conditions (iv) with following condition, we obtain the notion an extremal region with constraints on the right.

  (iv)$' \ S_i = \{z\in S_{i-1} \ | \ h_{g_i}(z) \geq h_i-q_i , z \ \mathrm{is \ connected \ to } \ I_i \ \mathrm{by \ a \  zigzag \ from} \ \mathcal{R}_{right}\}$, for all $1\leq i\leq n$,
  
  \end{defn}

  Notice that conditions (i)-(v) imply that for all $1\leq i\leq n-1$, $I_{i+1}$ is connected to $I_i$ with a zigzag from $\mathcal{R}_{left}\}$. 
  
  \medskip
  
   In the definition of partial order above, by replacing $R(g_1, \ldots g_n; q_1, \ldots , q_n)$ with $R(g_1, \ldots g_n; q_1, \ldots , q_n; I_1, \dots , I_n)$ everywhere we define {\em a partial order with  constraint on the left (right)}. We will use partial orders both with and without constraints. Unless said otherwise ``a partial order" will have no constraint. 
   
  \medskip

  We now fix the sequence $(g_1, \ldots , g_L)$ and study suitable short segments with respect to the order $\prec $ that it defines:
  
  \medskip
  
   \begin{defn} Let $I = I^k(x) = \{u, v, w\}\subset F_1$ be a suitable short segment where $u$ is the leftmost and $w$ is the rightmost element; and let $\prec $ be a linear order on $P(F_1)$. We will say $I$ is of 
   
   type 1, if $v \prec u$ and $v\prec w$;
   
   type 2, if $u \prec v$ and $w\prec v$;
   
   type 3, if $u \prec v\prec w$;
   
   type 4, if $w\prec v\prec u$.
   
   If $T$ is an elementary piece with root at $I = I^k(x)$, we say $T$ is of type 1 (or type 2, 3, 4), if $I$ is of type 1 (or type 2, 3, 4 respectively). 
   \end{defn}
   
   \medskip
   
   We also need the following notions
   
   \medskip
   
   \begin{defn} [exceptional regions] Let $I = \{u_1, \ldots , u_m\} \subseteq I(x)\subset P(F_1)$ where elements of the segment are listed from leftmost to rightmost. We say $I$ is positively exceptional if  $u_1\prec u_2\prec \ldots \prec u_m$; and negatively exceptional if $u_m\prec u_{m-1}\prec \ldots \prec u_1$; we say $I$ is exceptional if it is either positively or negatively exceptional. Similarly, we say a region $R\subseteq F_1$ is exceptional, if either all segments $I\subseteq P(R)$ are positively exceptional or they are all negatively exceptional.
 \end{defn}  
     
  \medskip
  
  \begin{defn} [neighbor segments] Let $I^{k'}(x'), I^{k''}(x'')$ be suitable short segments in $F_1$. We say they are neighbors if there exists an unsuitable segment $I^k(x)\subseteq F_1$ and $u, v\in I^k(x)$ such that $u\xi ^p\in I^{k'}(x'), v\xi ^p\in I^{k''}(x'')$. $I^k(x)$ will be called {\em a connecting segment}.     
  \end{defn}
  
  \medskip
  
  \begin{rem} In other words, two short intervals are neighbors if there is a quasi-zigzag of length 4 connecting them. Notice that by condition $D$-(iii) the connecting segment is unique.
  \end{rem}
  
  \medskip
  
  \begin{defn} [connected components] A sub-region of $F_1$ is called connected if any two suitable short segments in it are connected with a quasi-zigzag respecting the tiling  $\{I^k(x)\}_{x\in X, 0\leq k\leq 359}$. If $R$ is a maximal connected sub-region of $F_1$ then $R\sqcup N(R)$ is called a connected component of $F_1$.  
  \end{defn}
  
  \bigskip
  
   Now we are ready to start the proof of the theorem. Without loss of generality we may and will assume that $F_1$ is connected (otherwise we consider each connected component separately). Then, by Remark 2.3 if $u, v\in F_1$ belong to the suitable intervals of $F_1$ then either both $h(u), h(v)$ are even or both are odd; without loss of generality again we may and will assume that $h(u)$ is even whenever $u$ belongs to a suitable segment of $F_1$.    
   
   \medskip
   
   We can choose a sequence $g_1, \ldots , g_L$ such that the partial order $\prec $ imposed on $P(F_1)$ is linear and the partial order (still denoted with $\prec $) imposed on $X(F_1)$ is strongly linear, i.e. given any two suitable short intervals indexed by elements of $X(F_1)$, one is strongly bigger than the other. Then, without loss of generality and by shifting the tiling of the group $\Gamma $ by $k\in \{1, 2\}$ units to the right if necessar$X(F_1)$ y, we may assume that either $F_1$ is an exceptional region or there exists a finite collection $\{I^{k_i}(x_i) \ | \ 1\leq i\leq M\}$ of suitable segments such that

	\medskip
	
	 (i) $I^{k_i}(x_i)\subset F_1$ for all $1\leq i\leq M$;
	
	\medskip
	
	 (ii) $I^{k_1}\prec \ldots \prec I^{k_M}$, i.e. $(k_1, x_1)\prec \ldots \prec (k_M, x_M)$
	 
	 \medskip
	 
	 (iii) $M \geq \frac{1}{100|B_{50p}|}|F_1|$;
	 
	 \medskip
	 
	 and one of the following conditions hold:
	 
	 \medskip
	 
	 {\em Case A:} all short segments $I^{k_i}(x_i)$ are of type 1;
	 
	\medskip
	 
	 {\em Case B:} all short segments $I^{k_i}(x_i)$ are of type 3, and all of them have a neighbor of different type which precedes the short segment $I^{k_i}(x_i)$.
	
	\medskip

	{\em Case C:} all short segments $I^{k_i}(x_i)$ are of type 4 and all of them have a neighbor of different type which precedes the short segment $I^{k_i}(x_i)$.
			
	\bigskip
	
	We will first assume that for all short intervals $I^{k_i}(x_i), 1\leq i\leq M$ one of the cases $A, B, C$ holds.
	
	\medskip
	
	We will start the proof by describing the base (more precisely, the first step) of the inductive process:

	\medskip
	
	Let $R_1 = R_0(g_1, \ldots g_{N})$. Since the order $\prec $ is linear, $R_1$ contains only one short segment. Let $I^{k_0}(x_0)$ be this short segment. Then $I^{k_0}(x_0)$ has an absolutely minimal element, and let $z$ be this element.

	\medskip
	
	Then we build a complete labeled $\eta $-normal special elementary piece $T_1$ with $root(T_1) = I^{k_0}(x_0)$, and the starting element at $N(z)$ such that $T_1$ is successfully related to $I^{k_0}(x_0)$ and let $Y_1' = \{(x,k)\in Y_1 \ |  \ I^k(x)\cap S(T_1) \neq \emptyset \},  F_2 = \displaystyle \mathop{\sqcup }_{(x,k)\in Y_0\backslash Y_1'}I^k(x)$, and carry the inductive process by applying it to $F_2$. 	
	
 \medskip

 We continue the process of building special labeled $\eta $-normal pieces $T_1, T_2, \ldots , T_N$ and regions $F_1, F_2, \ldots , F_N$ inductively such that 

 \medskip

 (i) the sequence $T_1, T_2, \ldots , T_N$ is complete;

 \medskip

 (ii) $Y_i' = \{(x,k)\in Y_0 \ | \ I^k(x)\cap \displaystyle \mathop{\sqcup }_{1\leq j\leq i-1} S(T_j) \neq \emptyset \}$ for all $1\leq i\leq N$;

 \medskip

(iii) $F_i = \displaystyle \mathop{\sqcup }_{(x,k)\in Y_0\backslash Y_i'}I^k(x)$, for all $1\leq i\leq N$;   

 \medskip
 
 (iv) for every $i\in \{1, \ldots , N\}, T_i$ is a special elementary piece such that $root(T_i)$ is the biggest short interval $I$ in the set $P(F_i)$ with respect to the linear order $\prec $, and $T_i$ is successfully related to $I$ relative to $\displaystyle \mathop{\sqcup }_{1\leq j\leq i-1}S(T_j)$ (so, in particular, the starting element of $T_i$ is $N(z)$ where $z$ is the biggest element of $I$);   

 \medskip
 
 (v) for all $1\leq i\leq N$, $S(T_i)$ is a subset of an extremal region $R_0(g_1, \ldots , g_L)$ of $F_i$.

\medskip 

 (vi) $F_1\subseteq \displaystyle \mathop{\sqcup }_{1\leq i\leq N}S(T_i)$. 

 \bigskip

  Now, for every $i\in \{1, \ldots , M\}$, let $N_i\in \{1, \ldots , N\}$ be such that $I^{k_i}(x_i)$ is the biggest element of $P(F_{N_i})$ (notice that $N_i$ is uniquely determined and $N_1 < N_2 < \ldots < N_M$).
  
  \medskip
  
  Let also $I^{k_i}(x_i) = (u_i, v_i, w_i), 1\leq i\leq M$ where $u_i$ is the leftmost and $w_i$ is the rightmost element. 
  
  \medskip
  
  Now, we {\bf \em assume Case A} by the assumption of which, we have $v_i \prec u_i, v_i \prec w_i$. Let $N(v_i)\in I^k(x)$ for some $(k,x)\in Y_1$.  Let $I^k(x) = \{\alpha , N(v_i), \beta \}$ where $\alpha , \beta $ the leftmost and rightmost elements respectively. 
  
  \medskip
   
  Then $\alpha $ is successfully related to $u_i$ and $\beta $ is successfully related to $w_i$. Then $\alpha , \beta \in \displaystyle \mathop{\sqcup }_{1\leq j\leq N_i-1}S(T_j)$ because otherwise $I^{k_i}(x_i)$ is not the biggest element of $P(F_{N_i})$. But then, because of completeness, both of $\alpha $ and $\beta $ are starting elements of some elementary pieces $T_{j_1}, T_{j_2}$ where $j_1 < N_i, j_2 < N_i$. Hence, the special piece $T_{N_i}$ is successful relative to $\displaystyle \mathop{\sqcup }_{1\leq j\leq N_i-1}S(T_j)$, moreover, the sequence $T_1, \ldots , T_{N_i}$ is successful at $I^k(x)$.

  \medskip
  
  Now, let us {\bf \em assume Case B:} (Case $C$ is similar to Case $B$). Then $u_i\prec v_i\prec w_i, 1\leq i\leq M$. Let $I' = (u, v, w)$ be the short suitable segment preceding $I^{k_i}$ which is a neighbor of $I^{k_i}$ and has a different type. Let us assume it has type 4. (other cases are similar/easier). 
  
  \medskip
  
  Let also $u$ be the leftmost and $w$ be the rightmost element of $I'$, and $I_1, I_2, I_3$ be the unsuitable short segments containing $N(u), N(v), N(w)$ respectively. Then one of these segments is the connecting segment of $I^{k_i}$ and $I'$. 
  
  \medskip
  
  Let $I_3$ be a connecting segment. Since $I'$ has type 4, both elements in $I_3\backslash \{N(w)\}$ are successfully related to $u$ hence $I'$ cannot precede $I^{k_i}$ which contradicts our assumption.
  
  \medskip
  
  Let now $I_2$ be the connecting segment. Then $I_2 = (\alpha , N(v), \beta )$ where $\alpha , \beta $ are the the leftmost and the rightmost elements respectively. Then $\alpha $ is successfully related to $u$ and $\beta = N(w_i)$, hence $T_{N_i}$ is successful relative to $\displaystyle \mathop{\sqcup }_{1\leq j\leq N_i-1}S(T_j)$, and $T_1, \ldots , T_{N_i}$ is successful at $I_2$.

 \medskip
 
 Finally, let $I_1$ be the connecting segment. Then $N(u)\in I_1$ is the starting element of one of the $T_j, 1\leq j < N_i$, moreover, the set $(I_1\backslash \{N(u)\})\cap \displaystyle \mathop{\sqcup }_{1\leq j\leq N_i-1}S(T_j)$ either contains the starting element of one of the $T_{N_i}$ or it consists of a pair vertex of $T_{N_i}$. Hence again $T_{N_i}$ is successful relative to $\displaystyle \mathop{\sqcup }_{1\leq j\leq N_i-1}S(T_j)$, and $T_1, \ldots , T_{N_i}$ is successful at $I_1$.
  
 \medskip
 
 Then, by Proposition \ref{prop:+-} we obtain that  $|\partial _KF| > \frac{1}{900|B_{50p}|}|F|$ which is a contradiction.
 
 \medskip 
 
  Thus we obtain that either all but at most $M$ of the suitable short segments of $F_1$ are of type 3 or all but at most $M$ of the suitable short segments of $F_1$ are of type 4. But notice that if $\Gamma $ is amenable then for all $\epsilon > 0$, it admits $(F, \epsilon )$-Følner sets as well (i.e. one can replace $K$ with $F$). This implies the following intermediate proposition which is interesting in its own right
  
  \begin{prop} \label{prop:exceptional} If $\Gamma $ is amenable and satisfies conditions (A) and (D), then for all $\epsilon > 0$, $\Gamma $ admits $(K, \epsilon)$-Følner sets $F$ which is either positively exceptional or negatively exceptional.
  \end{prop}

 \bigskip

   By Proposition \ref{prop:exceptional}  we may assume that there exists a sequence $C$ which induces a linear order $\prec _1$ on $F_1$ which is either positively exceptional or negatively exceptional. On the other hand, by conditions $D$-(iii) and $(D)$-(iv), there exists a sequence $C'$ which induces a linear order $\prec '$ on $X(F_1)$ and $P(F_1)$ with right constraint such that no suitable segment $I^k(x)$ in $F_1$ is of type 3, and there exists a sequence $C''$ which induces a linear order on $X(F_1)$ and $P(F_1)$ with left constraint such that no suitable segment $I^k(x)$ in $F_1$ is of type 4. Thus we have one of the following two cases:

   \medskip
   
   Case 1. There exist sequences $C_1, C_2$ inducing linear orders $\prec _1, \prec _2$ on $X(F_1)$ and $P(F_1)$ such that $\prec _1$ is positively oriented, and $\prec _2$ is a linear order with right constraint.
   
   \medskip
   
    Case 2. There exist sequences $C_1, C_2$ inducing linear orders $\prec _1, \prec _2$ on $X(F_1)$ and $P(F_1)$ such that $\prec _1$ is negatively oriented and $\prec _2$ is a linear order with left constraint.
   
   \medskip
   
   These two cases are symmetric and we will be assuming we are in Case 1. Then all suitable short segments in $F_1$ are of type 3 with respect to the order $\prec _1$, and no suitable short segment in $F_1$ is of type 3 with respect to the order $\prec _2$. Because of the right constraint, then, any suitable short segment $I$ in $F_1$ is either of type 2 or of type 4, and the rightmost element of $I$ is the least element of it. 
   
   \medskip

   If we have at least $\frac{1}{900|B_{50p}|}|F|$ distinct pairs $(I_{2j-1}, I_{2j}), 1\leq j\leq m$ of suitable short segments in $F_1$ such that these pairs are connected with mutually non-interfering zigzags then again we obtain a complete sequence of labeled $\eta $-normal quasi-GBTs covering $F_1$ and being successful in at least $m$ non-suitable short segments. Then by Proposition \ref{prop:+-} we again obtain a contradiction. Thus we may assume that there exists $\Delta \subset Y_0$ where $\frac{|\Delta |}{|Y_0|} > 1- \frac{1}{100|B_{50p}|}$ such that, with respect to the ordering $\prec _2$,  either for all $(x,k)\in Y_1$ the suitable segment $I^k(x)$ is of type 2 or for all $(x,k)\in Y_1$ the suitable segment $I^k(x)$ is of type 4. Let us assume the latter case (the former case is very similar).
   
   \medskip
   
    Now, we will be working with both of the orderings $\prec _1$ and $\prec _2$. Let $g_1$ be the least element of $C_1$, $Y_{odd}^{(1)} = \{(x,k)\in Y_0 : \forall u\in I^k(x), h_{g_1}(u)\in 2\Z _{odd}\}$ and $Y_{even}^{(1)} = \{(x,k)\in Y_0 : \forall u\in I^k(x), h_{g_1}(u)\in 2\Z _{even}\}$. Then either $\frac{|Y_{odd}^{(1)}|}{|Y_0|} \geq \frac{1}{2}$ or $\frac{|Y_{even}^{(1)}|}{|Y_0|} \geq \frac{1}{2}$. Without loss of generality we may assume that $\frac{|Y_{odd}^{(1)}|}{|Y_0|} \geq \frac{1}{2}$. Then we let $Y_1 = Y_{odd}^{(1)}$.
    
    \medskip
    
    Let $$G_1 = \displaystyle \mathop{\cup } _{(x,k)\in Y_{odd}^{(1)}}(I^k(x)\cup N(I^k(x))), H_1 =  \displaystyle \mathop {\cup }_{(x,k)\in Y_{even}^{(1)}}(I^k(x)\cup N(I^k(x))).$$
    
    Let us observe that if for some $(x,k)\in Y_{odd}^{(1)}$, the suitable segment $I^k(x)$ contains an element $u$ with $h_{g_1}(u\xi ^{-p}) < h_{g_1}(u)$ then the following three conditions hold: 
    
    (i) $u$ is the rightmost element of $I^k(x)$, 
   
    (ii) if for some $(x',k')\in Y_0$ the suitable segment $I^{k'}(x')$ is a neighbor of $I^k(x)$ and $I^k(x) \prec _1 I^{k'}(x')$, then $I^{k'}(x')\subseteq H_1$, moreover, $h_{g_1}(w') = h_{g_1}(w)+2p$ for some (hence for all)  $w\in I^{k}(x), w'\in I^{k'}(x')$,
    
    (iii) if for some $(x',k')\in Y_0$ the suitable segment $I^{k'}(x')$ is a neighbor of $I^k(x)$, contains an element $v$ with $h_{g_1}(v\xi ^{-p}) < h_{g_1}(v)$ and $I^{k'}(x')\prec _1 I^k(x)$, then $I^{k'}(x')\subseteq H_1$, moreover, $h_{g_1}(w') = h_{g_1}(w)-2p$ for some (hence for all)  $w\in I^{k}(x), w'\in I^{k'}(x')$,
    
    \medskip
    
    For all non-negative integers $n$, let also $$G_{1,n} = \{g\in G_1 : h_{g_1}(g) = 4n-2\}\sqcup N(\{g\in G_1 : h_{g_1}(g) = 4n-2\}),$$ \ $$H_{1,n} = \{g\in H_1 : h_{g_1}(g) = 4n\}\sqcup N(\{g\in H_1 : h_{g_1}(g) = 4n\}).$$ 
    
    \medskip

    Then for all distinct $i, j\in \Z_{+}$, 
    there is no segment with non-trivial intersections with both $G_{1,i}$ and $G_{1,j}$, and similarly, there is no segment with non-trivial intersections with both $H_{1,i}$ and $H_{1,j}$; in addition, if $i-j\notin \{0, 1\}$, then again no segment intersects both $G_{1,i}$ and $H_{1,j}$. However, if $i-j\in \{0, 1\}$, $u\in G_{1,i}, v\in H_{1,j}$ and $u, v$ belong to the same short segment $I^k(x)$, then this segment is necessarily unsuitable; moreover, if $i-j = 0$, then $v$ is the rightmost element of $I^k(x)$, but if $i-j = 1$ then $v$ is not the rightmost element of $I^k(x)$.   
    
    \bigskip

    Let $$\Omega _1 = \{(x,k)\in Y_{odd}^{(1)} : \mathrm{there \  exists} \ u \in I^k(x) \ \mathrm{ such \ that} \ h_{g_1}(u\xi ^{-p}) < h(u)\}.$$  
    
    Let $T_1^{(1)}, \dots , T_{n_1}^{(1)}$ be the sequence of special elementary pieces of $H_1$ induced by the linear order $\prec _2$. This sequence forms a complete sequence of labeled $\eta $-normal elementary GBTs. Let also $m_1 = |\Omega _1|$ and $T_{n_1+1}^{(1)}, \dots , T_{n_1+m_1}^{(1)}$ be all elementary pieces successfully related to elements of $\Omega _1$. 
    
    \medskip

    Let $n_1+1 \leq i\leq n_1 + m_1$ and $v: = start (T_i)\in I^k(x), (x,k)\in X$. Then $v\in G_{1,n}$ for some $n$ and $I^k(x)\backslash \{v\}\subseteq H_{1,n-1}$ hence $I^k(x)\backslash \{v\}\subseteq \displaystyle \mathop{\sqcup }_{1\leq k\leq n_1}S(T_k)$. Let $j$ be the smallest number such that $(I^k(x)\backslash \{v\})\cap  S(T_j) \neq \emptyset $. Then there exists a suitable short segment $I^l(z)$ forming a vertex $w$ of $T_j$ (i.e. $S(w) = I^l(z)$) such that for a vertex $w'\in n(w)$ we have $S(w')\cap ((I^k(x)\backslash \{v\}) \neq \emptyset $. We make a key observation that if $(l,z)\in \Delta $ then $S(w') = (I^k(x)\backslash \{v\})$ (in particular,  $w'$ is a pair vertex) and the central element of $w'$ is the leftmost element of the segment $I^k(x)$.  
    
    \medskip
    
    Now lets us consider the other vertices of $T_i$; let $v_1, v_2 \in \Gamma \backslash \{v\}$ forming end vertices of $T_i$ where $v_1$ is the leftmost element of some unsuitable segment $I^{k'}(x')$ and $v_2$ is the middle element of some unsuitable segment $I^{k''}(x'')$. Let $j_1, j_2$ be the smallest numbers such that $(I^{k'}(x')\backslash \{v_1\})\cap  S(T_{j_1}) \neq \emptyset $ and $(I^{k''}(x'')\backslash \{v_1\})\cap  S(T_{j_2}) \neq \emptyset $. Then there exists suitable short segments $I^{l_1}(z_1), I^{l_2}(z_2)$ forming vertices $w_1, w_2$ of $T_j$ such that for some vertices $w_1'\in n(w_1), w_2'\in n(w_2)$ we have $S(w_1')\cap ((I^{k'}(x')\backslash \{v_1\}) \neq \emptyset $ and $S(w_2')\cap ((I^{k'}(x')\backslash \{v_2\}) \neq \emptyset $. Then we again make a a key observation that if $(l_1,z_1)\in \Delta $ then $S(w_1') = (I^{k'}(x')\backslash \{v_1\})$ (in particular,  $w_1'$ is a pair vertex) and the central element of $w_1'$ is the rightmost element of the segment $I^{k'}(x')$. On the other hand, if   $(l_2,z_2)\in \Delta $ then $S(w_2')$ is an end vertex consisting either the rightmost or the leftmost element of $I^{k''}(x'')$. 
    
    \medskip 
    
    Let now $C_1 = (g_1, g_2, g_3, \dots )$ where $g_i$ is the $i$-th least element of $C_1$ for all $i\geq 1$. We now define $Y_2 = Y_1\backslash \Omega _1$ and inductively define sequences $G_i, H_i, \Omega _i, Y_{i}, 1\leq i\leq r$ and the special elementary pieces $T_1^{(i)}, \dots , T_{n_i+m_i}^{(i)}$ and stop the process when for all suitable short segments $I^k(x), (x,k)\in Y_0$ there exists a ball $B_{u}(10p)$ of radius $10p$ centered at some $u\in I^k(x)$ such that $B_{u}(10p)\cap \displaystyle \mathop{\cup }_{1\leq i\leq r}H_i \neq \emptyset $ (we will have $r \leq |C_1|$). Then    
    
    \medskip 
    
    Now, let $\overline{Y} = Y_1\backslash (\displaystyle \mathop{\sqcup }_{1\leq k\leq r}\Omega _k)$. Let also $$F_{tail} = (\displaystyle \mathop{\sqcup }_{(x,k)\in \overline{Y}}I^k(x)\cup N(I^k(x)))\sqcup (F_1\backslash \displaystyle \mathop{\cup }_{1\leq i\leq r}H_i).$$

 On the set $\overline{Y}$ we have a linear order $\prec _1$. This order induces a complete sequence of elementary special GBTs $\Theta _1, \dots , \Theta _s$ which cover $F_{tail}$. 
    
    \medskip
    
    The quasi-GBTs $T_1^{(i)}, \dots , T_{n_i}^{(i)},  T_{n_i}^{(i)+1}, \dots , T_{n_i+m_i}^{(i)}, 1\leq i\leq r, \Theta _1, \dots , \Theta _s$ form a complete sequence. Among these, the special elementary pieces $T_{n_i}^{(i)+1}, \dots , T_{n_i+m_i}^{(i)}, 1\leq i\leq r$ are successfully related to elements of $\Omega _i$. Besides, we already observed that for any of these elementary pieces, if it is rooted at a suitable short segment $I_1 = \{z_1, z_2, z_3\}$ with $z_1$ being the leftmost and $z_3$ the rightmost element and a special elementary piece from the list $T_1^{(i)}, \dots , T_{n_i}^{(i)}$ is rooted at neighboring suitable short segments $I_2 = I^k(x)$ with $(k,x)\in  \Delta $ the the following holds true: if $I$ is the unsuitable segment connecting $I_1$ and $I_2$ (such a segment is unique) then the sequence $T_1^{(i)}, \dots , T_{n_i}^{(i)},  T_{n_i}^{(i)+1}, \dots , T_{n_i+m_i}^{(i)}, 1\leq i\leq r, \Theta _1, \dots , \Theta _s$ is successful at $I$ in both of the cases when $N(z_1)\in I$ and $N(z_3)\in I$ (in case $N(z_2)\in I$ the sequence can even be unsuccessful). Thus we obtain that there exists $m$ unsuitable short intervals $I_1, \dots , I_m$ and $n$ unsuitable intervals $I_{m+1}, \dots , I_{m+n}$ in $F_1$ with $m-n \geq \frac {1}{100|B_{50p}|}|F|$ such that $T_1^{(i)}, \dots , T_{n_i}^{(i)},  T_{n_i}^{(i)+1}, \dots , T_{n_i+m_i}^{(i)}, 1\leq i\leq r, \Theta _1, \dots , \Theta _s$ is successful at all $I_k, 1\leq k\leq m$, unsuccessful at all  $I_k, m+1\leq k\leq m+n$ and neither successful nor unsuccessful at all other unsuitable short intervals of $F_1$. Then by Proposition \ref{prop:+-},  $|\partial _KF| \geq \frac {|F_1|}{9|B_{70p}|}$. Contradiction.
    
   \bigskip

  \medskip
  
  Notice that if none of the Case A, Case B, Case C can be guaranteed, then it means that, loosely speaking, most of $F_1$ (100 percent in the limit) consists of regions which is either positively exceptional or negatively exceptional. Moreover, the places where positively and negatively exceptional regions meet have insignificant cardinality (zero percent in the limit). This is a rather extreme case, however, it seems difficult (if possible) to take care of it with just the inequalities of condition $(D)$ if we exclude the conditions $h(g^{(q)}\xi ^p) = h(g^{(q)}+p, 1\leq q\leq 2$. On the other hand, having the group $\mathbf {F}$ in mind as a major application, the linear order $\prec $ seems to be in agreement with a bi-order of $\mathbf{F}$ and this seems to cooperate with the possibility of this extreme case. Since  $\mathbf {F}$ does not have many interesting quotients, it is impossible to achieve one of the cases A-C, by taking quotient; it is clear that, for example, the following type condition rules out the possibility described above (i.e. the existence of large enough exceptional regions), so it guarantees the existence of one of the cases A, B, C:
  
  \medskip
  
  $(\ast )$ for all $\epsilon \in \{-1,1\}$, there exist $k\in \mathbb{N}, r\in \Z ,$ and $m_1, \ldots , m_k\in \N\backslash \{1,2\}, n_1, \ldots , n_k\in \{1, 2\}$ such that $$\eta ^{\epsilon m_1}\xi ^{-p}\eta ^{-\epsilon n_1}\xi^{p}\ldots \eta ^{\epsilon m_k}\xi ^{-p}\eta ^{-\epsilon n_k}\xi^{p} = \eta ^r.$$
  
  \medskip
  
  Indeed, by condition $(\ast )$, there cannot be an exceptional sub-region $R$ of $F_1$ which contains a ball of radius $R = 2pk + (m_1 + n_1) + \ldots + (m_k + n_k) + r$ (notice that we have two such quantities, one for each value of $\epsilon $; the radius $R$ can be taken as the maximum of these quantities). 
  
  Hence we obtain a result which is interesting in itself:
  
  \medskip
  
  \begin{thm} \label{prop:extra} If $\Gamma $ satisfies conditions (A), $(\ast )$ and the following weaker version of condition $(D)$:
  
  there exists an odd integer $p\in \Zodd $ such that for all $g\in \Gamma $, 

  (i) for at least one $\delta \in \{0,1\} $, the equality
      $h(g\eta^{\delta }\xi ^{-p}) = h(g)+ p \ (1)$ is satisfied.
\medskip

  (ii) if for some $\delta _0\in \{0,1\}$ the equality $(1)$ is not
satisfied then $h(g\eta^{\delta _0}\xi ^{p}) = h(g)+ p$, and for all $i\in \{1, 2\}, j\in [1,i], \epsilon \in \{-1, 1\}$, the inequality $h(g\eta^{\delta _0-\epsilon i}\xi ^{-p}\eta ^{\epsilon j}\xi ^p) \geq h(g\eta^{\delta _0-\epsilon i}\xi ^{-p}) + p$ holds.

 (iii) for all $m, n \in \{-2, -1, 1, 2\}$ the equality $h(\xi ^{-p}\eta ^{m}\xi ^{p}\eta ^{n}\xi ^{-p})= 3p$ holds.
  
  Then $\Gamma $ is non-amenable. $\square $
  \end{thm}
   
 \begin{rem} Notice that in the above statement of Theorem \ref{prop:extra} we have somewhat weakened condition D-(ii). The additional strength of this claim is needed only in our arguments beyond Proposition \ref{prop:exceptional}.  
 \end{rem}

\newpage
	
\begin{center} {\Large Part 2: Application to R.Thompson's group F} \end{center}

\vspace{1cm}

In 1965 Richard Thompson introduced a remarkable infinite group $\F$
that has two standard presentations: a finite presentation with two
generators and two relations, and an infinite presentation that is
more symmetric.
\[ \F \cong \langle A, B \mid [AB^{-1}, A^{-1}BA] = 1, [AB^{-1},
A^{-2}BA^2] = 1 \rangle\]
\[ \cong \langle X_0, X_1, X_2, \ldots \mid X_n X_m = X_m X_{n+1},\
\forall\ n>m \rangle\]

\noindent
Basic properties of $\F$ can be found in \cite{BS} and \cite{CFP}. Non-amenability of $\F $ is conjectured by Ross Geoghegan in \cite{Ge}. The standard isomorphism between the two presentations of $\F$ identifies $A$ with
$X_0$ and $B$ with $X_1$.  For convenience, let $\A = \{A,B\} = \{X_0,
X_1\}$, let $\X =\{X_0, X_1, X_2, \ldots \}$ and let $\free(\A)$ and
$\free(\X)$ denote the free groups of rank 2 and of countably infinite
rank with bases $\A$ and $\X$, respectively.

\vspace{1cm}

\section{Normal Form}

Recall the following basic fact about elements in free groups.

\begin{prop}[Syllable normal form]\label{prop:syllable}
  There is a natural one-to-one correspondence between non-trivial
  elements in $\free(\X)$ and words $X_{n_1}^{e_1} X_{n_2}^{e_2}\cdots
  X_{n_k}^{e_k}$ with nonzero integer exponents and distinct adjacent
  subscripts.
\end{prop}

A word of the form described in Proposition \ref{prop:syllable} is
called the \emph{syllable normal form} of the corresponding element in
$\free(\X)$.  The terminology refers to the language metaphor under
which an element of $\X$ is a \emph{letter}, a finite string of
letters and their formal inverses is a \emph{word}, and a maximal
sub-word of the form $X_n^e$ is a \emph{syllable}.  

 \medskip

 We will use the following result on normal forms for elements in Thompson's group
$\F$.  For a proof of this result see \cite{CFP}.

\begin{thm}[Thompson normal forms]\label{thm:thompson}
  Every element in $\F$ can be represented by a word $W$ in the form
  $X_{n_0}^{e_0} X_{n_1}^{e_1} \ldots X_{n_k}^{e_k} X_{m_l}^{-f_l}
  \ldots X_{m_0}^{-f_0}$ where the $e$'s and $f$'s are positive
  integers and the $n$'s and $m$'s are non-negative integers
  satisfying $n_0 < n_1 < \ldots < n_k$ and $m_0 < m_1 < \ldots <
  m_l$.  If we assume, in addition, that whenever both $X_n$ and
  $X_n^{-1}$ occur, so does either $X_{n+1}$ or $X_{n+1}^{-1}$, then
  this form is unique and called the \emph{Thompson normal form} of
  this element.
\end{thm}

In order to cleanly describe the rewriting process used to convert an
arbitrary reduced word into its equivalent Thompson normal form (and
to explain the reason for the final restrictions), it is useful to
introduce some additional terminology.

\begin{defn}[Shift map]
  Let $S$ denote the map that systematically increments subscripts by
  one.  For example, if $W = X_{n_1}^{e_1}X_{n_2}^{e_2}\ldots
  X_{n_k}^{e_k}$, then $S(W)$ is the word $X_{n_1+1}^{e_1}
  X_{n_2+1}^{e_2} \ldots X_{n_k+1}^{e_k}$.  More generally, for each
  $i\in \N$, let $S^i(W)$ denote $i$ applications of the shift map to
  $W$. Thus, $S^i(W) = X_{n_1+i}^{e_1} X_{n_2+i}^{e_2} \ldots
  X_{n_k+i}^{e_k}$.  Note that this process can also be reversed, a
  process we call \emph{down shifting}, so long as all of the
  resulting subscripts remain non-negative.  Also note, that a shift
  of an odd word, up or down, remains an odd word.
\end{defn}

\begin{rem}[Rewriting words]
  Using the shift notation, the defining relations for $\F$ can be
  rewritten as follows: for all $n>m$, $X_n X_m = X_m S(X_n)$.  More
  generally, let $W = X_{n_1}^{e_1}X_{n_2}^{e_2}\ldots X_{n_k}^{e_k}$
  be a reduced word and let $\min(W)$ denote the smallest subscript
  that occurs in $W$, i.e. $\min(W) = \min\{n_1, n_2, \ldots, n_k\}$.
  It is easy to show that for all words $W$ with $\min(W) > m$, $W X_m
  = X_m S(W)$.  Similarly, $X_m^{-1} W = S(W) X_m^{-1}$. After
  iteration, for all positive integers $e$, $W X_m^e = X_m^e S^e(W)$
  and $X_m^{-e} W = S^e(W) X_m^{-e}$.  The reason for the extra
  restriction in the statement of Theorem \ref{thm:thompson} should
  now be clear.  If $W$ is in Thompson normal form with $\min(W)>m$,
  then $X_m S(W) X_m^{-1}$ is an equivalent word that satisfies all
  the conditions except the final restriction.
\end{rem}

\begin{defn}[Core of a word]\label{def:core}
  Let $W$ be a reduced word with $\min(W)=m$. By highlighting those
  syllables that achieve this minimum, $W$ can be viewed as having the
  following form:
  \[W = W_0 X_m^{e_1} W_1 X_m^{e_2} \ldots X_m^{e_l} W_l\]
  where the $e$'s are nonzero integers, each word $W_i$ is a reduced
  word with $\min(W_i) > m$, always allowing for the possibility that
  the first and last words, $W_0$ and $W_k$, might be the empty word.

  We begin the process of converting $W$ into its Thompson normal form
  by using the rewriting rules described above to shift each syllable
  $X_m^{e_i}$ with $e_i$ positive to the extreme left and each such
  syllable with $e_i$ negative to the extreme right.  This can always
  be done at the cost of increasing the subscripts in the subwords
  $W_i$.  If we let $pos$ and $neg$ denote the sum of the positive and
  negative $e$'s, respectively, then $W$ is equivalent in $\F$ to a
  word of the form $W' = X_m^{pos} W_0' W_1' \ldots W_l' X_m^{neg}$
  with $W_i'$ is an appropriate upward shift of the word $W_i$.  The
  appropriate shift in this case is the sum of the positive $X_m$
  exponents in $W$ to the right of $W_i$ plus the absolute value of
  the sum of the negative $X_m$ exponents in $W$ to the left of $W_i$.
  The resulting word $W_0' W_1' \ldots W_l'$ between $X_m^{pos}$ and
  $X_m^{neg}$ is called the \emph{core} of $W$ and denoted $\core(W)$.
\end{defn}

The construction of the core of a word, is at the heart of the process
that produces the Thompson normal form.

\begin{rem}[Producing the Thompson normal form]
  Let $W$ be a reduced word and let $W' = X_m^{pos} \core(W)
  X_m^{neg}$ be the word representing the same element of $\F$
  produced by the process described above.  If the first letter of
  $W'$ is $X_m$, the last letter is $X_m^{-1}$ and $\min(\core(W)) >
  m+1$ then we can cancel an $X_m$ and an $X_m^{-1}$ and downshift
  $\core(W)$ to produced an equivalent word whose core has a smaller
  minimal subscript.  We can repeat this process until the extra
  condition required by the normal form is satisfied with respect to
  the subscript $m$.  At this stage we repeat this entire process on
  the new core, the down-shifted $\core(W)$.  After a finite number of
  iterations, the end result is an equivalent word in Thompson normal
  form.
\end{rem}

 From the description of the rewriting process, the following
proposition should be obvious.

\begin{prop}[Increasing subscripts]\label{prop:increase}
  If $W$ is word with $\min(W)=n$ and a non-trivial Thompson normal
  form $W'$, then $\min(W')$ is at least $n$.  In particular, when
  $\min(W) >m$, the words $W$ and $X_m^e$, $e$ nonzero, represent
  distinct elements of $\F$.
\end{prop}

\begin{exmp}
  Consider the following word:
\[W = (X_2^{-3} X_5^2) X_0^4 (X_1^5 X_3^{-2}) X_0^{-1}
(X_1^7) X_0^{2} (X_3 X_4)\]
\noindent
It has $\min(W) = 0$, $pos = 6$, $neg = -1$.  Pulling the syllables
with minimal subscripts to the front and back produces the equivalent
word:
\[ W' = X_0^6 (X_8^{-3} X_{11}^2) (X_3^5 X_5^{-2}) (X_4^7) (X_4 X_5) X_0^{-1}\]
\noindent
with $\core(W) = X_8^{-3} X_{11}^2 X_3^5 X_5^{-2} X_4^8 X_5$.  Note
that we needed to combine two syllables in order for the core to be in
syllable normal form.  The process of reducing this to Thompson normal
form would further cancel an initial $X_0$ with a terminal $X_0^{-1}$
and down shift the core because $\min(\core(W))=3 > 1+1$.  The new
word is:
\[ W' = X_0^5 (X_7^{-3} X_{10}^2 X_2^5 X_4^{-2} X_3^8 X_4) \]
\noindent
and the new core is $X_7^{-3} X_{10}^2 X_2^5 X_4^{-2} X_3^8 X_4$.
\end{exmp}

\medskip

Now we are close to claim that $\F$ satisfies the conditions of Theorem \ref{prop:main}, but for that, first, we need to introduce the height function

\begin{defn} Let $w = X_{n_0}^{e_0} X_{n_1}^{e_1} \ldots X_{n_k}^{e_k} X_{m_l}^{-f_l}\ldots X_{m_0}^{-f_0}$ be a Thompson normal form of $W\in \F $.
Then we let 

 $h(W) = (e_1 + \ldots + e_k) + (f_1 + \ldots + f_l)$, if $n_0 = 0, m_0 = 0$;
 
 $h(W) = (e_0 + e_1 + \ldots + e_k) + (f_1 + \ldots + f_l)$, if $n_0 > 0, m_0 = 0$;
 
 $h(W) = (e_1 + \ldots + e_k) + (f_0 + f_1 + \ldots + f_l)$, if $n_0 = 0, m_0 > 0$;
 
 $h(W) = (e_0 + e_1 + \ldots + e_k) + (f_0 + f_1 + \ldots + f_l)$, if $n_0 > 0, m_0 > 0$;
 
 $h(W) = 0$, if $W = 1\in F$.
\end{defn}

  \bigskip

\begin{prop} \label{prop:verify} With $\eta = X_0, \xi = X_1, p = 1$, the group $\F $ satisfies property $(D)$.
\end{prop}

     \bigskip

{\bf Proof.} It is clear that the function $h: \F \rightarrow \N\cup \{0\}$ is subadditive, and $h(g) = h(g^{-1})$ for all $g\in \F $.   

\medskip

 We will verify the conditions (i), (iii), and condition (ii) for $\epsilon = 1$; for $\epsilon = -1$ it is verified similarly. We will also verify condition (iv), but for the values +1 and -1 of $\epsilon $, our arguments will be somewhat different from each other.

  \medskip
  
  Notice also that condition $(D)$-(ii) unites the following three claims: the equalities $h(g\eta^{\delta _0}\xi ^{p}) = h(g)+ p \ (\star _1)$ and $h(g\eta^{\delta _0}\xi ^{-p}) = h(g) - p \ (\star _2)$ and the inequality $h(g\eta^{\delta _0-\epsilon i}\xi ^{-p}\eta ^{\epsilon j}\xi ^p) \geq h(g\eta^{\delta _0-\epsilon i}\xi ^{-p}) + p \  (\star _3)$.
   Condition $(D)$-(iii) claims the inequality $h(\xi ^{-p}\eta ^{m}\xi ^{p}\eta ^{n}\xi ^{-p})\geq 3p \ (\star _4)$ 
     
     \medskip

   We will first concentrate on condition $(D)$-(i) and on the (in)equalities $(\star _1), (\star _2)$ and $(\star _3)$. 
  
  \medskip
  
 Let $g\in \F , W =  X_{n_0}^{e_0} X_{n_1}^{e_1} \ldots X_{n_k}^{e_k} X_{m_l}^{-f_l}\ldots X_{m_0}^{-f_0}$ be the Thompson normal form of $g$. For any $\epsilon, \delta \in \Z $, let $W(\epsilon , \delta ) = g\eta ^{\delta }\xi ^{\epsilon }$. 

\medskip

 If $h(g\eta ^{\delta }\xi ^{-1}) = h(g)+1$ for all $\delta \in \{0, 1\}$ then we have nothing to prove. So let $\delta _0$ be the biggest number in the set $\{0, 1\}$ such that $h(W(-p, \delta _0)) \leq h(g) \ (\dagger )$. Let also $\delta = \delta _0 - i, i\in \{1,2\}$ and $W =  X_{n_0}^{e_0} X_{n_1}^{e_1} \ldots X_{n_k}^{e_k} X_{m_l}^{-f_l}\ldots X_{m_0}^{-f_0}$ be the Thompson normal form of $g\eta ^{\delta _0} $. 
 
 \medskip
 
 We will consider several cases:

 \medskip
 
 {\em Case 1: $m_0 > 0$.} 
 
 \medskip
 
 Then, necessarily $n_0 = 1, m_0\geq 3, n_1\geq 3$. Then $$g\eta ^{\delta _0}\eta ^{-1}\xi ^{-p} = WX_0^{-1}X_1^{-1} = X_{n_0}^{e_0} X_{n_1}^{e_1} \ldots X_{n_k}^{e_k} X_{m_l}^{-f_l}\ldots X_{m_0}^{-f_0}X_2^{-1}X_0^{-1}$$ and the latter expression is in the normal form. Hence $h(g\eta ^{\delta _0}\eta ^{-1}\xi ^{-p}) = h(g) + p$, and this proves the claim $D$-(i). For condition $D$-(ii), we consider the following sub-cases:
 
 \medskip 
 
 a) For $i = 1, j = 1$, we have $$W\eta ^{-1}\xi ^{-1}\eta \xi = WX_0^{-1}X_1^{_1}X_0X_1 = WX_2^{-1}X_1 = X_1S(WX_2^{-1}).$$ Hence $$W\eta ^{-1}\xi ^{-1}\eta \xi = X_1X_{n_0+1}^{e_0} X_{n_1+1}^{e_1} \ldots X_{n_k+1}^{e_k} X_{m_l+1}^{-f_l}\ldots X_{m_0+1}^{-f_0}X_3^{-1}$$ and since $n_0 = 1$ the latter expression is in the normal form therefore we obtain that $$h(W\eta ^{-1}\xi ^{-p}\eta \xi ^p = h(g) + 2p = h(W\eta ^{-1}\xi ^{-p}) + p.$$ 
 
 \medskip 
 
 b) For $i = 2, j = 1$, we have $$W\eta ^{-2}\xi ^{-1}\eta \xi = WX_0^{-2}X_1^{-1}X_0X_1 = WX_3^{-1}X_0^{-1}X_1 = WX_3^{-1}X_2X_0^{-1}.$$ But $$WX_3^{-1}X_2X_0^{-1} = X_2S(WX_3^{-1})X_0^{-1} = X_2S(W)X_4^{-1}X_0^{-1}$$ and $$W\eta ^{-2}\xi ^{-1} = WX_3^{-1}X_0^{-2}$$ hence we obtain that $h(W\eta ^{-2}\xi ^{-p}\eta \xi ^p = h(W\eta ^{-1}\xi ^{-p}) + p$.

 \medskip 
 
 c) For $i=2, j = 2$, we have $$W\eta ^{-2}\xi ^{-1}\eta ^2\xi = WX_0^{-2}X_1^{-1}X_0^2X_1 = WX_3^{-1}X_1$$ and  $$W\eta ^{-2}\xi ^{-1} = WX_3^{-1}X_0^{-2}$$ therefore, again, $h(W\eta ^{-2}\xi ^{-p}\eta \xi ^p = h(W\eta ^{-1}\xi ^{-p}) + p$. 
 
 \bigskip 
 
 {\em Case 2: $m_0 = 0$.}
 
 \medskip

  In this case since $h(gX_0^{\delta _0}X_1^{-1}) = h(g) - 1$, we have that $gX_0^{\delta _0}X_1^{-1} = WX_1^{-1} = X_{n_0}^{e_0} X_{n_1}^{e_1} \ldots X_{n_s}^{e_s}\dots X_{n_k}^{e_k} X_{m_l}^{-f_l}\ldots X_{m_{r+1}}^{-f_{m_r+1}} X_{n_s}^{-1} X_{m_{r}}^{-f_{m_r}} \dots X_{m_0}^{-f_0}$ where $n_s > m_r$. moreover, if $X_{m_{r+1}}$ occurs (in the normal form $W$), then $m_{r+1}\geq n_s+2$, and similarly, if $X_{n_{s+1}}$ occurs, then $n_{s+1} \geq n_s + 2$. Then the verification of condition D-(ii) is done similarly as in Case 1.

  \medskip
  
  {\em Case 3: $W = X_{n_0}^{e_0} X_{n_1}^{e_1} \ldots X_{n_k}^{e_k}$, i.e. the negative part of the Thompson normal form is absent.}
  
  \medskip
  
   This case is not different from the previous cases except we observe that the inequality $h(g\eta ^{\delta }\xi ) = h(g) + 1$ holds for all $\delta \in \Z$.
	
	\medskip
	
 Since $p=1, \eta = X_0, \xi = X_1$, the inequality $(\star _2)$ simply follows from the fact that for all $g\in \F$ either $h(gX_1) = h(g) -1 $ or $h(gX_1) = h(g) + 1$. 
	
\medskip 

	For the inequality $(\star _4)$ of condition $D$-(iv), we notice that $$\xi ^{-p}\eta ^{i_1}\xi ^{p}\eta ^{i_2}\xi ^{-p}\eta ^{i_3}\xi ^{p} = X_0^rX_{n_1}^{-p}X_{n_2}^{p}X_{n_3}^{-p}X_{n_4}^{p}X_0^{s}$$ and the claim can be seen easily by a direct check.

	\bigskip 
	
	Now we will verify condition (iv). We will treat the cases $\epsilon = 1$ and $\epsilon = -1$ separately, and our arguments for these two cases will be different.  
	
	\medskip
	
	{\em Case \ $\epsilon = -1$.} 
	
	\medskip

	   Let $$u = \eta ^{j_0}\xi ^{-p}\eta ^{i_1}\xi ^{p}\eta ^{j_1}\ldots \xi ^{-p}\eta ^{i_k}\xi ^{p}\eta ^{j_k}$$ where $s_0 = j_0\in \{0,1,2\}$ and for all $1\leq q\leq k$ we have $r_q, s_q\in \{0, 1, 2\}, r_q\neq 2, s_k = 2$.  
	   
	   \medskip
	   
	   Let also $Q = \{q : 0\leq q\leq k, s_q = 2\} = \{q_1, q_2, \dots , q_r\}$ where $q_1 < q_2 < \dots < q_r$. Notice that since $s_k = 2$, we have $q_r = k$. Then we let $u_s = \xi ^{-p}\eta ^{i_1}\xi ^{p}\eta ^{j_1}, 1\leq s\leq k, u_0 = \eta ^{j_0}$, and $V_{s} = u_{q_{s}+1}\dots u_{q_{s+1}}, 1\leq s\leq r-1$. We also let $V_0 = u_0u_1\dots u_{q_1}$ (if $j_0 = 2$, then we let $V_0 = \eta ^2$). 
	   
	   \medskip 
	   
	   Since $\eta = X_0, \xi = X_1$, for each $s\in \{1, \dots , r\}$ we have either 
	   
	   \medskip 
	   
	   $$V_s = (X_1^{-1}X_0^{-1}X_1X_0^{-1})(X_1^{-1}X_0X_1X_0^{-1})^nX_0^2$$ or $$V_s = (X_1^{-1}X_0^{-2}X_1X_0)(X_1^{-1}X_0^{-1}X_1X_0)^nX_0.$$ 
	   
	   \medskip
	   
	    Then $$V_s = X_1^{-1}X_2(X_3^{-1}X_2)^n = X_3(X_4^{-1}X_3)^nX_1^{-1}$$  or $$V_s = X_1^{-1}X_3(X_2^{-1}X_3)^n = X_4(X_3^{-1}X_4)^nX_1^{-1}.$$ 
	   
	   \medskip
	   
	   On the other hand, $V_0 = X_0^eU$ where $e\in \{1,2\}$ and either $U = 1$ or $\min(U) \geq 2$. Then the normal form of $u$ equals $$W = X_{n_0}^{e_0} X_{n_1}^{e_1} \ldots X_{n_k}^{e_k} X_{m_l}^{-f_l}\ldots X_{m_0}^{-f_0}$$ where $m_0 = 1$ (and $f_0 > 0$), moreover, if $j_0\in \{1,2\}$ then $n_0 = 0, e_0 = j_0, n_1 > 1$, but if $j_0 = 0$, then $n_0 > 1$. Then $u\xi ^{-p}$ has a normal form  $$X_{n_0}^{e_0} X_{n_1}^{e_1} \ldots X_{n_k}^{e_k} X_{m_l}^{-f_l}\ldots X_{m_0}^{-f_0-1}$$ hence $h(u\xi ^{-p}) = h(u)+p$.

	\bigskip

	{\em Case \ $\epsilon = 1$.} 
	
	\medskip 
	
	 Let again $$u = \eta ^{j_0}\xi ^{-p}\eta ^{i_1}\xi ^{p}\eta ^{j_1}\ldots \xi ^{-p}\eta ^{i_k}\xi ^{p}\eta ^{j_k}$$ where $s_0 = j_0\in \{0,1,2\}$ and for all $1\leq q\leq k$ we have $r_q, s_q\in \{0, 1, 2\}, r_q\neq 0, s_k = 0$. Let also $$Q = \{q : 1\leq q\leq k, s_q = 0\} = \{q_1, q_2, \dots , q_r\}$$
	 
	 where $q_1 < q_2 < \dots < q_r$. Then we let again $u_s = \xi ^{-p}\eta ^{i_1}\xi ^{p}\eta ^{j_1}, 1\leq s\leq k, u_0 = \eta ^{j_0}$, and $V_{s} = u_{q_{s}+1}\dots u_{q_{s+1}}, 1\leq s\leq r-1$. We also let $V_0 = u_0u_1\dots u_{q_1}$ (if $j_0 = 0$, then we let $V_0 = 1$). 
	 
	 \medskip
	 
	 Notice that $u = V_0V_1\dots V_{r-1}$. Since $\xi = X_1, \eta = X_0, p = 1$, if $(u\xi ^{-p}) < h(u)+p$ then $h(V_0V_1\dots V_{r-2}\xi ^{-p}) < h(V_0V_1\dots V_{r-2})+p$ and inductively, we obtain that $h(V_0\xi ^{-p}) < h(V_0) + p$. Then $V_0 \neq 1$ hence $j_0 \neq 0$. Then either $$V_0 = (X_0X_1^{-1}X_0X_1)(X_0^{-1}X_1^{-1}X_0X_1)^{n-1}X_0^{-2}$$
	 or $$V_0 = (X_0^2X_1^{-1}X_0^{-1}X_1)(X_0X_1^{-1}X_0^{-1}X_1)^{n-1}X_0^{-1}$$ for some $n\geq 1$.
	 
	 After simplifying, in the former case we obtain $$V_0 = X_0^2(X_2^{-1}X_1)^nX_0^{-2} = X_0^2X_1^nX_{n+2}^{-1}X_{n+1}^{-1}\dots X_3^{-1}X_0^{-2}$$ whereas in the latter case we obtain that $$V_0 = X_0^2(X_1^{-1}X_2)^nX_0^{-2} = X_0^2X_3X_4\dots X_{n+2}X_1^{-n}X_0^{-2}.$$ Thus in both cases we have $h(V_0\xi ^{-p}) = h(V_0)+p$ which contradicts our assumption.  $\square $

\bigskip

\begin{thm} $\F $ satisfies all conditions of Theorem \ref{prop:main} therefore it is non-amenable.
\end{thm}

\medskip

{\bf Proof.} It is well known that $\F /[\F, \F]$ is isomorphic to $\Z ^2$. (See \cite{CFP}). We choose $\eta = X_0, \xi = X_1$. Then condition $(D)$ follows from Proposition \ref{prop:verify}


\vspace{1cm}

\bigskip 


\begin{thebibliography}{99}

\bibitem{BS} Brin, M, Squier, C. \ Groups of piecewise linear homeomorphisms of the real line.  \ {\em Inventiones Mathematicae} {\bf 79} (1985), no.3, 485-498.

\medskip 

\bibitem{CFP} Cannon, J.W, Floyd,W.J, Parry,W.R. \ Introductory notes on Richard Thompson's groups. \ {\em Enseign.  Math.}  (2)  {\bf 42}  (1996), \ no 3-4.

\medskip

\bibitem{F} Følner, E. \ On groups with full Banach mean value, {\em Math. Scand.} (1955) \ vol.3, 243-254.

\medskip 

\bibitem{Ge} Open problems in infinite-dimensional topology. Edited by Ross Geoghegan. The Proceedings of the 1979 Topology Conference. {\em Topology Proc.} {\bf 4} (1979), no.1, 287-338. 

\end{thebibliography}
\end{document}